\documentclass[12pt]{article}
\usepackage{amsmath,amsfonts,amssymb,amsthm,eucal}
\usepackage[dvips]{graphics}

\newcommand{\C}{\mathbb{C}}
\newcommand{\R}{\mathbb{R}}
\newcommand{\PP}{\mathbb{P}}

\newcommand{\M}{\mathcal{M}}
\newcommand{\Mb}{\overline{\mathcal{M}}}
\newcommand{\Mgn}{\M_{g,n}}
\newcommand{\Mgnb}{\Mb_{g,n}}

\newcommand{\bG}{\mathsf{G}}
\newcommand{\Ggnb}{\bG^{\ge 3}_{g,n}}
\newcommand{\bH}{\mathsf{H}}

\newcommand{\lan}{\left\langle}
\newcommand{\ran}{\right\rangle}

\DeclareMathOperator{\Spec}{Spec}

\DeclareMathOperator{\Aut}{Aut}
\DeclareMathOperator{\Hur}{Hur}
\DeclareMathOperator{\val}{val}
\DeclareMathOperator{\Prob}{Prob}
\DeclareMathOperator{\tr}{tr}

\begin{document}

\title{Random trees and moduli of curves}
\author{Andrei Okounkov\thanks{
 Department of Mathematics, University of California at
Berkeley, Evans Hall \#3840, 
Berkeley, CA 94720-3840. E-mail: okounkov@math.berkeley.edu}}
\date{}
\maketitle

\begin{abstract}
This is an expository account of the proof
of Kontsevich's combinatorial formula for
intersections on moduli spaces of curves
following the paper \cite{OP}. It is based
on the lectures I gave on the subject in 
St.\ Petersburg in July of 2001.  
\end{abstract}

These are notes from the lectures I gave in 
St.\ Petersburg in July of 2001. Our goal here is to give an 
informal introduction to intersection theory on
the moduli spaces of curves and its relation to 
random matrices and combinatorics. More
specifically, we want to explain the 
proof of Kontsevich's formula given in \cite{OP}
and how it is connected to other topics discussed 
at this summer school such as, for example, the 
combinatorics of increasing subsequences in 
a random permutations.  

These lectures were intended for an audience of mostly
analysts and combinatorialists interested in 
asymptotic representation 
theory and random matrices. This is very much reflected in both 
the selection of the material and its 
presentation. Since absolutely no background in geometry
was assumed, there is a long and very basic  discussion
of what moduli spaces of curves and intersection 
theory on them are about. We hope that a reader 
trained in analysis or combinatorics will get some
feeling for moduli of curves (without worrying too
much about the finer points of the theory, all but a few  of 
which were swept under the rug). 

Conversely, in the 
second, asymptotic, part of the text, I allowed myself
to operate more freely because the majority of the
audience was experienced in asymptotic analysis. 
Also, since many fundamental ideas such as e.g.\  the KdV equations for 
the double scaling limit of the Hermitian $1$-matrix
model were at length discussed in other lectures of
the school, their discussion here is much more 
brief that it would have been in a completely
self-contained course. A much more detailed treatment
of both geometry and asymptotics can be found in 
the paper \cite{OP}, on which my lectures were based. 

It is needless to say that, since this is an expository
text based on my joint work with Rahul Pandharipande, 
all the credit should be
divided while any blame is solely my responsibility.
Many people contributed
to the success of the St.~Petersburg summer school, but 
I want to especially thank A.~Vershik for organizing 
the school and for the invitation to
participate in it. I am grateful to NSF (grant DMS-0096246),
Sloan foundation, and Packard foundation for partial
financial support.

\section{Introduction to moduli of curves}

\subsection{}

Let me begin with an analogy. 
In the ideal world, the moduli spaces of curves
would be quite similar to the Grassmann varieties
$$
Gr_{k,n}=\{L\subset \C^n,\dim L=k \}
$$
of $k$-dimensional linear subspaces $L$ of an $n$-dimensional
space. While any such subspace $L$ is geometrically just
a $k$-dimensional vector space $\C^k$, nontrivial things can 
happen if $L$ is allowed to vary in families, and this
nontriviality is captured by the geometry of $Gr_{k,n}$.
 
A convenient formalization of the notion of a family of 
linear spaces parameterized by points of some base space $B$ 
is a (locally trivial) vector bundle over $B$. There is a 
natural \emph{tautological} 
vector bundle over the Grassmannian $Gr_{k,n}$ 
itself, namely the space
$$
\mathcal{L}=\{(L,v), v\in L \subset \C^n\} 
$$
formed by pairs $(L,v)$, where $L$ is a $k$-dimensional 
subspace of $\C^n$ and $v$ is a vector in $L$. Forgetting
the vector $v$ gives a map $\mathcal{L}\to Gr_{k,n}$ whose
fiber over $L\in Gr_{k,n}$ is naturally identified with
$L$ itself. 

Given any space $B$ and a map 
$$
\phi: B \to Gr_{k,n}
$$
we can form the pull-back of $\mathcal{L}$
$$
\phi^* \mathcal{L} =\{(b,v), b \in B, v\in C^n, v\in \phi(b)\}
$$
which is a rank $k$ vector bundle over $B$. For a compact base
$B$, in the $n\to\infty$ limit this becomes a bijection between
(homotopy classes) of maps $B\to Gr_{k,\infty}$ and 
(isomorphism classes) of rank $k$ vector bundles on $B$.  

In particular, one can associate to a vector bundle $\phi^* \mathcal{L}$
its characteristic cohomology classes obtained by pulling back
the elements of $H^*(Gr_{k,n})$ via the map $\phi$. 
Intersections of these classes describe the enumerative 
geometry of the bundle $\phi^* \mathcal{L}$. It is thus 
especially important to understand the intersection theory 
on the space $Gr_{k,n}$ itself --- and this leads to a very 
beautiful classical combinatorics, in particular, 
Schur functions play a central role
(see for example \cite{Fu}, Chapter 14).

\subsection{}
One would like to have a similar theory with families of 
linear spaces replaced by families of curves of some
genus $g$. That is, given a family $F$ of, say, smooth genus
$g$ algebraic curves parameterized by some base 
$B$ we want to have a natural map $\phi: B\to \M_g$ that  
captures the essential information about the family $F$. 
Here $\M_g$ is the \emph{moduli space} 
of smooth curves of genus $g$, that is, the space of
isomorphism classes of smooth genus $g$ curves. At this point it may be 
useful to be a little naive about what we mean by a
family of curves etc., in this way we should be able
to understand the basic issues without too many 
technicalities getting in our way. Basically, 
a family $F$ of curves with base $B$ is a ``nice'' morphism 
$$
\pi : F \to B
$$
of algebraic varieties whose fibers $\pi^{-1}(b)$, $b\in B$,
are smooth complete genus $g$ curves.
We want the moduli space $\M_g$ and the induced map 
$\phi: B \to \M_g$ to be also algebraic.  

We will see that the first difficulty with the above
program is that, in general, the family $F$ will not
be a pull-back of any universal family over $\M_g$. 
To get a sense of why this is the case we can cheat
a little and consider the (normally forbidden) case $g=0$.
Up to isomorphism, there is only one curve of genus $0$,
namely the projective line $\PP^1$. Hence the map $\phi$
in this case can only be the trivial map to a point. 
There exist, however, highly nontrivial families with 
fibers isomorphic to $\PP^1$ or even $\C^1$ as we, in fact, already
saw above in the example of the tautological rank $1$ bundle 
over $Gr_{1,n} \cong \PP^{n-1}$.

The reason why there exist locally trivial yet 
globally nontrivial families with fiber $\PP^1$ is that
$\PP^1$ has a large automorphism group which one
can use to glue trivial pieces in a nontrivial way. 
Basically, the automorphisms are the principal issue 
behind the nonexistence of a universal family of
curves over $\M_g$. The situation becomes manageable, if
not entirely perfect, once one can get the automorphism 
group to be finite (which is automatic for 
smooth curves of genus $g>1$). A standard way to achieve
this is to consider curves with sufficiently many 
marked points on them ($\ge 3$ marked points for $g=0$ and $\ge 1$ for
$g=1$). Since curves with marked points arise very naturally
in many other geometric situations, the moduli spaces $\M_{g,n}$
of smooth genus $g$ curves with $n$ distinct marked points
should be considered on equal footing with the moduli
spaces $\M_g$ of plain curves.

\subsection{}\label{sM11}

As the first example, 
let us consider the space $\M_{1,1}$ of genus $g=1$ curves
$C$ with one parked point $p\in C$. By Riemann-Roch the 
space of $H^0(C,\mathcal{O}(2p))$ of meromorphic functions on $C$
with at most double pole at $p$ has dimension $2$. Hence, in 
addition to constants, there exists  a (unique up to 
linear combinations with constants) nonconstant meromorphic
function 
$$
f: C \to \PP^1\,,
$$
with a double pole at $p$ and no other poles (this 
is essentially the Weierstra\ss\ function 
$\wp{}$). Thus, $f$ defines at a $2$-fold branched 
covering of $\PP^1$ doubly ramified over $\infty \in\PP^1$.
For topological reasons, it has three additional 
ramification points which, after normalization, we
can take to be $0,1$ and some $\lambda\in \PP^1\setminus\{0,1,\infty\}$. 

\begin{figure}[!hbt]
\centering
\scalebox{.5}{\includegraphics{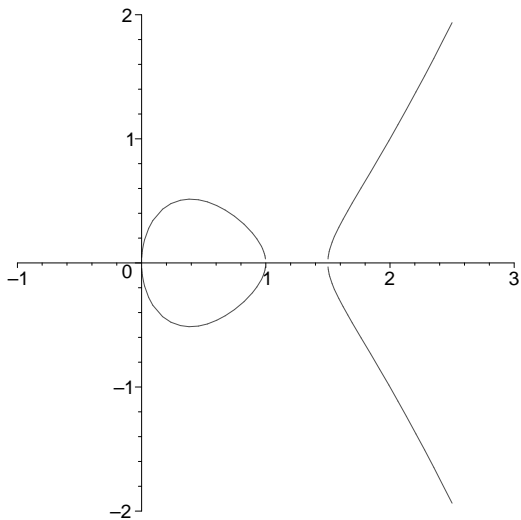}}
\caption{The $\lambda=\frac32$ member of the family \eqref{F3}.}
\label{fig1}
\end{figure}

Now it is easy to show that such a curve $C$ must be
isomorphic to the curve 
\begin{equation}
C\cong\{y^2= x(x-1)(x-\lambda)\}\in \PP^2\label{Cyx}
\end{equation}
in such a way that the point $p$ becomes the unique point at 
infinity and the function $f$ becomes the coordinate 
function $x$. It follows that
every smooth pointed $g=1$ curves occurs in the following 
family of curves
\begin{equation}
F=\{(x,y,\lambda), y^2= x(x-1)(x-\lambda)\} \label{F3}
\end{equation}
with the base
$$
B=\{\lambda\}= \PP^1\setminus\{0,1,\infty\}\,, 
$$
where the  marked point is the point $p$ at infinity. 
For example, the curve \eqref{Cyx} corresponding
to $\lambda=\frac32$ is plotted in Figure \ref{fig1}.

However, a given curve $C$ occurs in the family \eqref{F3} more 
than once. Indeed, we made arbitrary choices when we normalized 
the 3 critical values of $f$ to be $0$, $1$, and 
$\lambda$, respectively. At this stage, 
we can choose any of $6=3!$ possible assignments 
which makes the symmetric group $S(3)$ act on the
base $B$ preserving the isomorphism class of the fiber. 
Concretely, this group is generated by involutions 
$$
\lambda \mapsto 1-\lambda\,,
$$
which interchanges the roles of $0$ and $1$, and by 
$$
\lambda \mapsto 1/\lambda\,,
$$
exchanging the roles of $1$ and $\lambda$. It can be shown
that two  members of the family $F$ are isomorphic if and 
only if they belong to the same $S(3)$ orbit. Thus, the
structure map $\phi: B \to \M_{1,1}$ should be just 
the quotient map
$$
\phi: B \to B/S(3)= \Spec \C[B]^{S(3)} \,.
$$
Here $\C[B]^{S(3)}$ is the algebra of $S(3)$-invariant
regular functions on $B$. This algebra is a polynomial 
algebra with one generator, the traditional choice for
which is the following 
$$
j(\lambda) = 256\, \frac{(\lambda^2 - \lambda +1)^3}
{\lambda^2 (\lambda-1)^2} \,. 
$$
Thus, $\M_{1,1}$ is simply a line 
$$
\M_{1,1} = \Spec \C[j]  \cong \C  \,.
$$

It is now time to point out that the family $F$ is not a 
pull-back $\phi^*$ of some universal family over $\M_{1,1}$. 
The simplest way to see this is to observe that the
group $S(3)$ fails to act on $F$. Indeed,
let us try to lift the involution $\lambda\to 1-\lambda$
from $B$ to $F$. There are two ways to do this, namely
$$
(x,y,\lambda) \mapsto (1-x,\pm i y,1-\lambda)\,,
$$
neither of which is satisfactory because the square
of either map
$$
(x,y,\lambda) \mapsto (x,-y,\lambda)
$$
yields, instead of identity, a 
 nontrivial automorphism of every curve in 
the family $F$. One should also observe that 
both choices act by a nontrivial automorphism on the fiber 
over the fixed point $\lambda=\frac 12 \in B$. In fact,
the fibers of $F$ over fixed points of a transposition
and a 3-cycle in $S(3)$, respectively,
 (with $j(\lambda)=1728$ and $j(\lambda)=0$, 
resp.) are precisely the curves with extra
large automorphism groups (of order $4$ and $6$, resp.).

The existence of an nontrivial automorphism of
every pointed genus $1$ curve leads to the somewhat 
unpleasant necessity to consider every point of $\M_{1,1}$
as a ``half-point'' in some suitable sense in order to 
get correct enumerative predictions. Again, 
automorphisms make the real world not quite ideal.

While it is important to be aware of these
automorphism issues (for example, to understand
how intersection numbers on moduli spaces can be 
rational numbers), there is no need to be pessimistic
about them. In fact, by allowing spaces more general 
than algebraic varieties (called \emph{stacks}) 
one can live a life in which $\M_{g,n}$ is smooth and with a universal family
over it. This is, however, quite technical and will 
remain completely outside the scope of these lectures.

\subsection{}\label{Fr}

\begin{figure}[!hbt]
\centering
\scalebox{.5}{\includegraphics{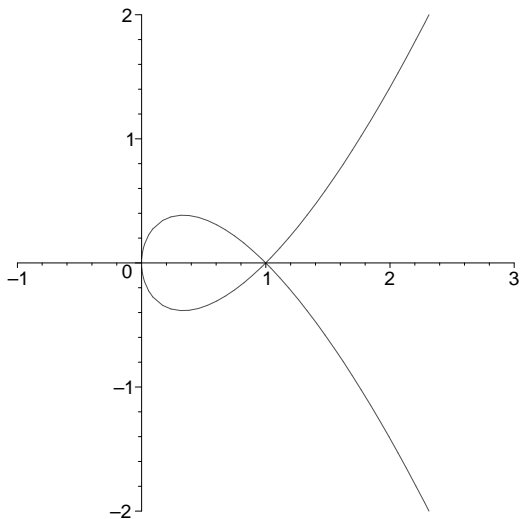}}
\caption{The nodal $\lambda=1$ curve in the family \eqref{F3}.}
\label{fig2}
\end{figure}

Clearly, the space $\M_{1,1}\cong \C$ is not compact. 
The $j$-invariant of the curve \eqref{Cyx} goes
to $\infty$ as the parameter $\lambda$ approaches the 
three excluded points $\{0,1,\infty\}$. As $\lambda$
approaches $0$ or $1$, the curve  $C$ acquires a nodal singularity; 
for example, for $\lambda=1$ we get the curve
$$
y^2 = x(x-1)^2 
$$
plotted in Figure \ref{fig2}.
It is natural to complete the 
family \eqref{F3} by adding the corresponding nodal cubics for 
$\lambda \in \{0,1\}$.  All plane cubic with a node
being isomorphic, the function $j$ extends
to a map 
$$
j: \C \to \PP^1/S(3)  = \Mb_{1,1} \cong \PP^1 
$$
to the moduli space of $\Mb_{1,1}$ of curves of 
arithmetic genus $1$ with at most one node and a smooth
marked point.\footnote{
For future reference we point out that something rather
different happens in the family \eqref{F3} as $\lambda\to\infty$.
Indeed, the equation 
$$
\lambda^{-1} \, y^2 = x(x-1)(\lambda^{-1} x - 1)
$$
becomes $x(x-1)=0$ in the $\lambda\to\infty$ limit,
which means that we get three lines (one of which is
the line at infinity), all three of them intersecting
in the marked point at infinity. In other words, the 
fiber of the family \eqref{F3} at $\lambda=\infty$ is very 
much not the kind of curve by which we want to 
compactify $\M_{1,1}$. This problem can be cured, but in 
a not completely trivial way, see below.}

In general, it is very desirable to have a nice 
compactification for the noncompact spaces $\M_{g,n}$. 
First of all, interesting families of curves over a complete base
$B$ are typically forced to have singular fibers over
some points in the base (as in the example above). 
Fortunately, as we will see below, 
it often happens that precisely these special fibers contain key 
information about the geometry of the family. 
Also, since eventually we will be interested in 
intersection theory on the moduli spaces of curves,
having a complete space can be a significant advantage. 

\begin{figure}[!hbt]
\centering
\scalebox{.3}{\includegraphics{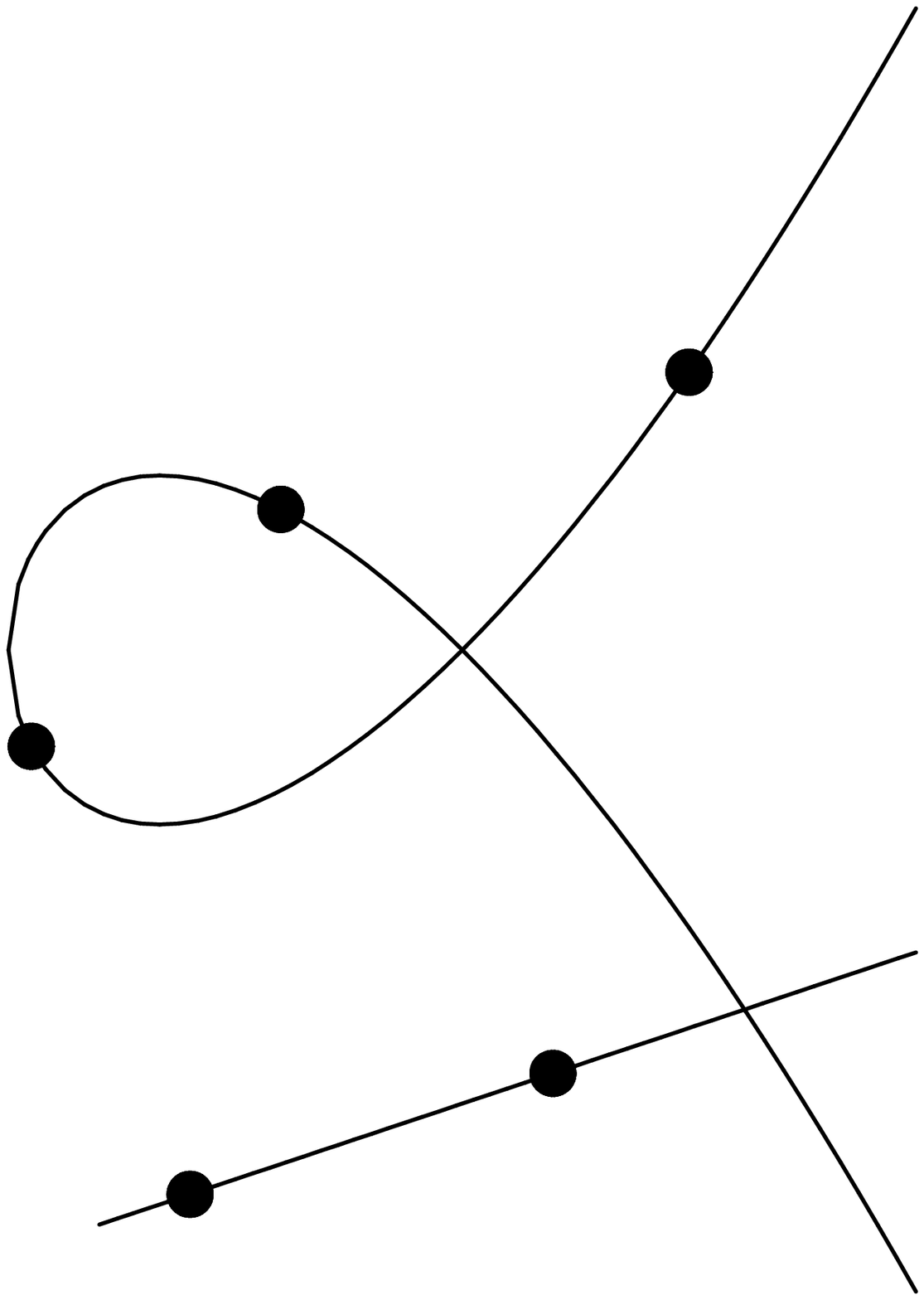}}
\caption{A boundary element of $\Mb_{1,5}$}
\label{fig3}
\end{figure}

A particularly remarkable compactification $\Mgnb$ of $\Mgn$
was constructed by Deligne and Mumford. The
space $\Mgnb$ is the moduli space of \emph{stable} curves
$C$ of arithmetic genus $g$ with $n$  distinct marked points. 
Stability, by definition, means  that 
the curve $C$ is complete and connected with
at worst nodal singularities, all marked points are 
smooth, and that $C$, together with marked points,
admits only finitely many automorphisms. In practice,
the last condition means that every rational 
component of the normalization of $C$ should have
at least 3 special (that is, lying over marked or 
singular points of $C$)
points. Observe that, in particular, the curve $C$ is allowed
to be reducible. A typical stable curve can be 
seen in Figure \ref{fig3}. 

\subsection{}\label{scomp}

Those who have not seen this before are probably left
wondering how it is possible for $\Mgnb$ to be compact.
What if, for example, one marked point $p_1$ on 
some fixed curve $C$ approaches 
another marked point $p_2\in C$ ?  
We should be able to assign some meaningful stable limit to such
a $1$-parametric family of curves, but it is somewhat
nontrivial to guess what it should be. 

A family of curves with a $1$-dimensional base $B$ is
a surface $F$ together with a map $\pi: F\to B$ whose
fibers are the curves of the family. Marking $n$ points
on the curves means giving $n$ sections of the map
$\pi$, that is, $n$ maps 
$$
p_1,\dots,p_n: B \to F\,,
$$ 
such that 
$$
\pi(p_k(b))=b\,,\quad b\in B\,,\quad  k=1,\dots,n\,.
$$
We will denote by 
$$
S_i=p_i(B)
$$
the trajectories of 
the marked points on $F$; they are curves on the 
surface $F$. 

\begin{figure}[!hbt]
\centering
\scalebox{.4}{\includegraphics{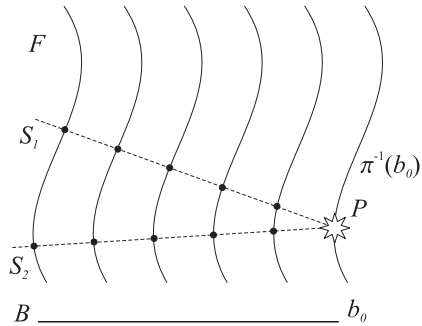}}
\caption{A family with colliding marked points}
\label{fig4}
\end{figure}

\begin{figure}[!hbt]
\centering
\scalebox{.4}{\includegraphics{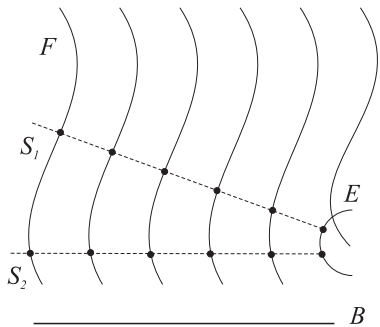}}
\caption{Same family, except with collision point blown up}
\label{fig5}
\end{figure}

Now suppose we have a $1$-dimensional family of $2$-pointed curves 
such that at some bad point $b_0$ of the base $B$ 
we have $p_1(b_0)=p_2(b_0)$, that
is, over this point two marked points hit each other, see 
Figure \ref{fig4}, and therefore the  fiber $\pi^{-1}(b_0)$ is not a 
stable $2$-pointed curve.
It is quite easy to repair this family: just blow up the
offending (but smooth) point $P=p_1(b_0)=p_2(b_0)$ 
on the surface $F$. Let 
$$
\sigma : \widetilde{F} \to F 
$$
be the blow-up at $P$. Then
$$
\widetilde{\pi} = \pi \circ \sigma: \widetilde{F} \to B
$$
is new family of curves with base $B$. Outside $b_0$ this 
the same family as before, whereas the fiber $\widetilde{\pi}^{-1}(b_0)$
is the old fiber $\pi^{-1}(b_0)$ plus the exceptional 
divisor $E=\sigma^{-1}(P)\cong \PP^1$ of the blow-up, see Figure \ref{fig5}.
Assuming the sections $S_1$ and $S_2$ met each other and 
the fiber $\pi^{-1}(b_0)$ transversally at $P$, as in Figure \ref{fig4},
the marked points
on $\widetilde{\pi}^{-1}(b_0)$ are two distinct point on the
exceptional divisor $E$. Therefore, $\widetilde{\pi}^{-1}(b_0)$
is a stable $2$-pointed curve which is the stable limit 
of the curves $\pi^{-1}(b)$ as $b\to b_0$. 

\begin{figure}[!hbt]
\centering
\scalebox{.5}{\includegraphics{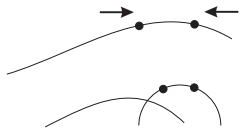}}
\caption{Bubbling off a projective line}
\label{fig6}
\end{figure}

To summarize, if one
marked point on a curve $C$ approaches another then $C$ bubbles
off a projective line $\PP^1$ with these two points on it as
in Figure \ref{fig6}. 

More generally, if $F$ is any family of curves with a smooth
1-dimensional base $B$ that are stable except over one offending point
$b_0\in B$ then after a sequence of blow-ups and blow-downs and,
possibly, after passing to a branched covering of $B$, one can 
always arrive at a new family with all fibers stable.  Moreover,
the fiber over $b_0$ in this family is determined uniquely. 
This process is 
called the \emph{stable reduction} and how it works is explained, for 
example,  in \cite{HM}, Chapter 3C. 

In particular, there exists a stable reduction of the family
\eqref{F3} which, as we saw in Section \ref{Fr}, fails to 
have a stable fiber over the point $\lambda=\infty$ in the
base. This is an example where only blow-ups and 
blow-downs will not suffice, that is, a base change is
necessary. 

\subsection{}

The topic of this lectures is intersection theory on 
the Deligne-Mumford spaces $\Mgnb$ and, specifically,
intersections of certain divisors $\psi_i$ which 
will be defined presently. It was conjectured by
Witten \cite{W} that a suitable generating function 
for these intersections is a $\tau$-function for the 
Korteweg--de Vries hierarchy of differential 
equations. This conjecture was motivated by an analogy
with matrix models of quantum gravity, where the 
same KdV hierarchy appears (this was already discussed
in other lectures at this school). The KdV 
equations were deduced by Kontsevich in \cite{K} from 
an explicit combinatorial formula for the 
intersections of the $\psi$-classes (see also, for example, 
\cite{D} for more information about the connection to 
the KdV equations). The main goal 
of these lectures is to explain a proof of this combinatorial 
formula of Kontsevich following the paper \cite{OP}.

The definition of the divisors $\psi_i$ is the 
following. A point in $\Mgnb$ is a stable curve
$C$ with marked points $p_1,\dots,p_n$. By
definition, all points $p_i$ are smooth points
of $C$, hence the tangent space $T_{p_i} C$ to 
$C$ at $p_i$ is a line. Similarly, we have 
the cotangent lines $T^*_{p_i} C$, $i=1,\dots,n$. As the point
$(C,p_1,\dots,p_n)\in \Mgnb$ varies, these cotangent lines
$T^*_{p_i} C$ form $n$ line 
bundles over $\Mgnb$. By definition, $\psi_i$
is the first Chern class of the line bundle
$T^*_{p_i} C$. In other words, it is the divisor of 
any nonzero section of the line bundle $T^*_{p_i} C$.

\subsection{}\label{ssi}

To get a better feeling for these classes let us 
intersect $\psi_i$ with a curve in $\Mgnb$. The 
answer to this question should be a number. Let $B$ be a curve.
A map $B\to\Mgnb$ is morally equivalent to a 
$1$-dimensional family of curves with base $B$ 
(in reality we may have to pass to 
a suitable branched covering of $B$ to get 
an honest family. \footnote{We already saw an 
example of this in Section \ref{sM11}. Indeed,
$\Mb_{1,1}$ is itself a line $\PP^1$. However,
in order to get an actual family over it
we have to go to a branched covering. As always, the 
automorphisms are to blame.}) So, let us
consider a family $\pi: F\to B$ of stable pointed curves
with base $B$ and the induced map
$\phi:B\to \Mgnb$. As usual, the marked points
$p_1,\dots,p_n$ are sections of $\pi$ and 
$$
S_i=p_i(B)\,, \quad i=1,\dots,n\,,
$$
are disjoint curves on the surface $F$. 
\begin{figure}[!hbt]
\centering
\scalebox{.4}{\includegraphics{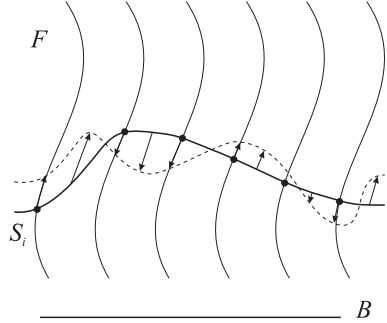}}
\caption{The self-intersection of $S_i=p_i(B)$ on $F$.}
\label{fig7}
\end{figure}
A section $s$ of $\phi^*(T_{p_i})$ is 
a vector field on the curve $S_i$ which is 
tangent to fibers of $\pi$ and, hence, $s$ is a section 
of the normal bundle to $S_i\subset F$,  see Figure \ref{fig7}. 
The degree of this normal bundle 
 is the self-intersection of the curve
$S_i$ on the surface $F$, that is, 
$$
\deg\,(s)=(S_i,S_i)_F \,,
$$ 
where $(s)$ is the divisor of $s$. In other words, 
$$
\int_{B} c_1\left(\phi^*(T_{p_i}) \right) = (S_i,S_i)_F \,,
$$
where $c_1$ denotes the $1$st Chern class. Dually, we have 
$$
\int_{\phi(B)} \psi_i  = - (S_i,S_i)_F\,. 
$$

We will now use this formula to compute the intersections
of $\psi_i$ with $\Mgnb$ in the cases when 
the space $\Mgnb$ is itself $1$-dimensional. Since
\begin{equation}
\dim\Mgnb = 3g -3 + n \label{dimMgnb}\,,
\end{equation}
this happens for $\Mb_{0,4}$ and $\Mb_{1,1}$.

\subsection{} \label{sM04}

The space $\Mb_{0,4}$ is easy to understand. After all,
there is only one smooth curve of genus $0$, namely
$\PP^1$. Moreover, any 3 distinct points of $\PP^1$ can
be taken to  the points $\{0,1,\infty\}$ by an 
automorphism of $\PP^1$ (in particular, this means that
$\Mb_{0,3}$ is a point). After we identified the first
three marked points with $\{0,1,\infty\}$, we can take
any point $x\in \PP^1 \setminus \{0,1,\infty\}$ as the
fourth marked point. Thus the locus of smooth curves
in $\Mb_{0,4}$ is isomorphic to $\PP^1 \setminus \{0,1,\infty\}$.
Singular curves are obtained as we let $x$ approach
the 3 excluded points $\{0,1,\infty\}$, which, by the 
process described in Section \ref{scomp}, bubbles off a new
$\PP^1$ with two marked points on it. This completes
the description of $\Mb_{0,4}\cong \PP^1$. 

In addition, this gives a description of the universal 
family over $\Mb_{0,4}$ (oh yes, in genus $0$ it does exist !). 
Take $\PP^1\times  \PP^1$ with coordinates $(x,y)$.
The map $(x,y)\to x$ with 4 sections 
$$
p_1(x)=(x,0)\,, \quad p_2(x)=(x,1)\,, \quad
p_3(x)=(x,\infty)\,, \quad p_4(x)=(x,x)\,,
$$
defines a family of 4-pointed smooth genus $0$ curves
for $x \in \PP^1 \setminus \{0,1,\infty\}$. 
\begin{figure}[!hbt]
\centering
\scalebox{.4}{\includegraphics{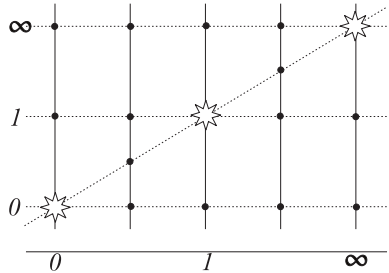}}
\caption{The trivial family $\PP^1\times \PP^1$ with sections
$p_1,\dots,p_4$.}
\label{fig8}
\end{figure}
The section $p_4$ collides with the other three at the 
points $(0,0)$, $(1,1)$, and $(\infty,\infty)$, see Figure \ref{fig8}.
To extend this family over all of $\PP^1$, we blow up these
collision points as in Section \ref{scomp} and get the surface $F$ shown in 
Figure \ref{fig9}. 
\begin{figure}[!hbt]
\centering
\scalebox{.4}{\includegraphics{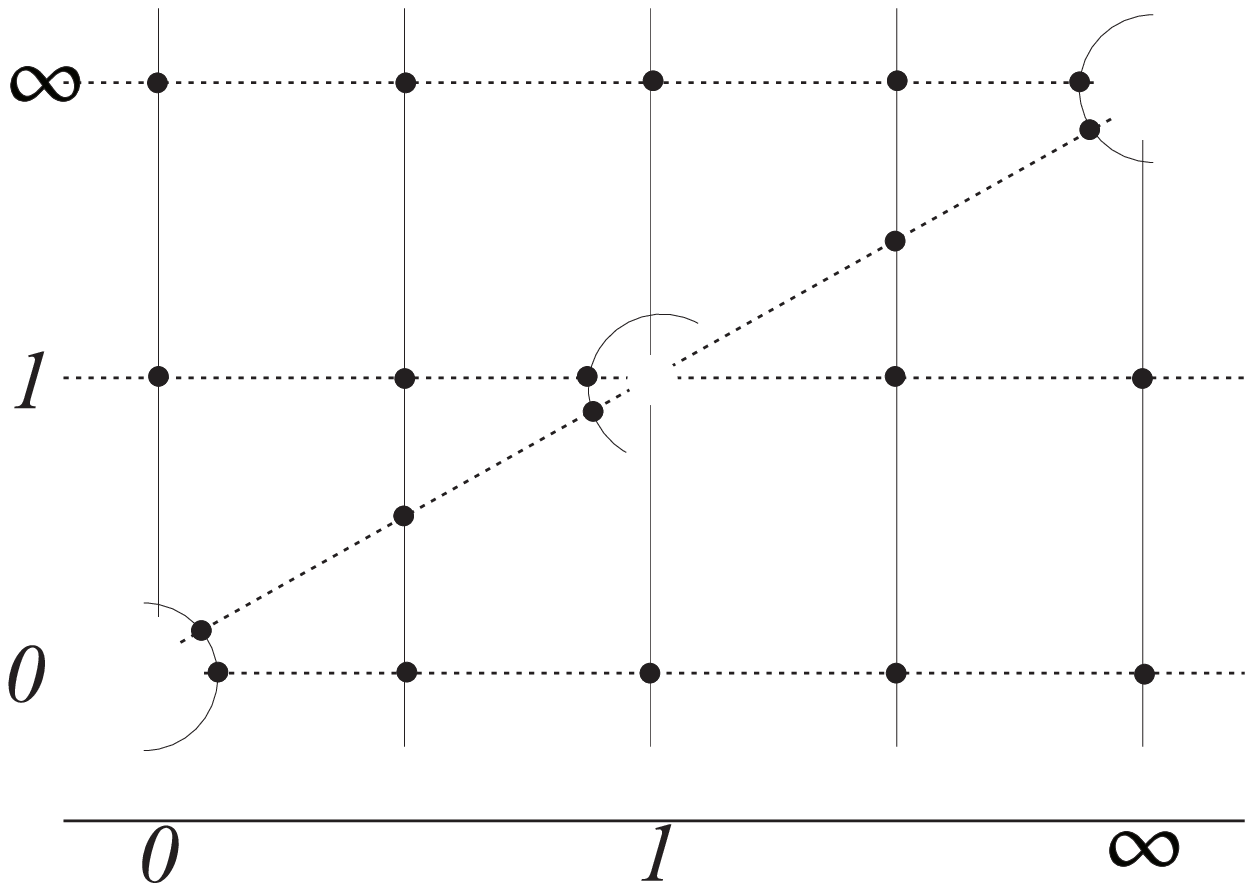}}
\caption{The universal family over $\Mb_{0,4}$.}
\label{fig9}
\end{figure}
The closures of the curves $p_1(x),\dots,p_4(x)$
give $4$ sections which are now disjoint everywhere.

Incidentally, this surface $F$ can be naturally identified
with $\Mb_{0,5}$ and, more generally, for any $n$ 
there exists a natural map 
$\Mb_{0,n+1}\to \Mb_{0,n}$ giving the universal
family over $\Mb_{0,n}$. This map forgets the 
$(n+1)$st marked point and, if the curve becomes
unstable after that, blows down all unstable components. 

Now let us compute $\int_{\Mb_{0,4}} \psi_1$ using the 
recipe given in Section \ref{ssi}. Recall that $S_1\subset F$ denotes
the closure of the curve $\{(x,0)\}$, $x\ne 0$ in $F$ (a.k.a.\ 
the proper transform of the corresponding 
curve in $\PP^1\times  \PP^1$). Let $E$ denote the preimage of
$(0,0)$ under the blow-up, that is, let $E$ be the exceptional
divisor. The self-intersection of $S_1$ with
any  curve $\{y=c\}$, $c\ne 0$, on 
$\PP^1\times  \PP^1$ is clearly zero. Letting $c\to 0$ we
get 
$$
(S_1,E+S_1)=0 \,.
$$
Since, obviously, $(S_1,E)=1$ we conclude that
$$
\int_{\Mb_{0,4}} \psi_1 = - (S_1,S_1)=-(-1)=1 \,.
$$

\subsection{}

\begin{figure}[!hbt]
\centering
\scalebox{.6}{\includegraphics{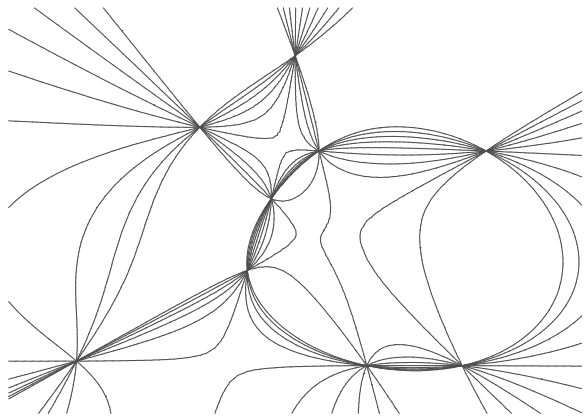}}
\caption{A pencil of plane cubics.}
\label{fig10}
\end{figure}

Now let us analyze the integral $\int_{\Mb_{1,1}} \psi_1$. 
In the absence of a universal family, we have to look for
another suitable family to compute this integral. A particularly
convenient family can be obtained in the following way. 
Consider the projective plane $\PP^2$ with affine coordinates
$(x,y)$. Pick two generic cubic polynomials $f(x,y)$ and 
$g(x,y)$ and consider the family of cubic curves
\begin{equation}
F=\{(x,y,t), f(x,y)- t \, g(x,y) = 0\} \subset \PP^2\times \PP^1\,,
\label{pencilF}
\end{equation}
with base $B=\PP^1$ parameterized by $t$. The cubic curves
$f(x,y)=0$ and $g(x,y)=0$ intersect in 9 points $p_1,\dots,p_9$
and we can choose any of those points as the marked point
in our family. An example of such family of plane cubics
is plotted in Figure \ref{fig10}. 

\begin{figure}[!hbt]
\centering
\scalebox{.6}[0.5]{\includegraphics{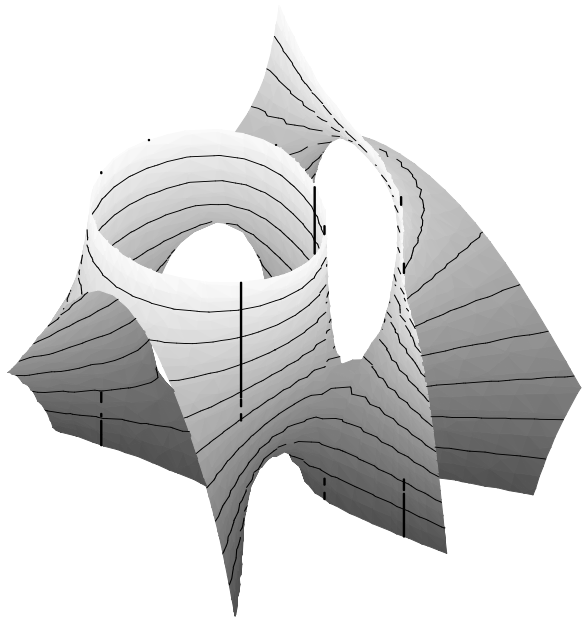}}
\caption{A fragment of the surface \eqref{pencilF}}
\label{fsur}
\end{figure}

Our first observation is that the surface $F$ is the blow-up
of $\PP^2$ at the points $p_1,\dots,p_9$ (which are
distinct for generic $f$ and $g$). Indeed, we have a rational
map
$$
\PP^2 \owns (x,y)\to \left(x,y,\frac{g}{f}\right) \in F \,,
$$
which is regular away from $p_1,\dots,p_9$. Each of the points
$p_i$ is a transverse intersection of $f=0$ and $g=0$, which is
another way of saying that at those points the differentials
$df$ and $dg$ are linearly independent. Thus, this map identifies
$F$ with the blow-up of $\PP^2$ at $p_1,\dots,p_9$. The graph of
the function $\frac{g(x,y)}{f(x,y)}$ is shown in Figure \ref{fsur}.
This graph goes vertically over the points $p_1,\dots,p_9$
that are blown up.

Since the section $S_1$ is exactly one of the exceptional 
divisors of this blow-up, arguing as in 
Section \ref{sM04} above we find that
$$
(S_1,S_1) = -1 \,.
$$
It does not mean, however, that we are done with the computation 
of the integral, because the induced map $\phi:B \to \Mb_{1,1}$ is very
far from being one-to-one. In fact, set-theoretically, the 
degree of the map $\phi$ is $12$ as we shall now see. 

To compute the degree of $\phi$ we need to know how many times
a fixed generic elliptic curve appears in the family $F$. This is a
classical
computation. First, one can show 
that the singular cubic is generic enough. Then we claim that, as $t$
varies, there will be precisely $12$ values of $t$ that produce
a singular curve. There are various ways to see this. For example,
the singularity of the curve is detected by vanishing of the 
discriminant. The discriminant of a cubic polynomial
is a polynomial of degree 12 in its coefficients, 
hence a polynomial of degree 12 in $t$.

A alternative way to obtain this number $12$ is to compute
the Euler characteristic of the surface $F$ in two different ways.
On the one hand, viewing $F$ as a blow-up, we get
$$
\chi(F) = \chi(\PP^2) + 9(\chi(\PP^1)-1) = 12  \,.
$$
On the other hand, $F$ is fibered over $B$ and the generic
fiber is a smooth elliptic curve whose Euler characteristic is 0.
The special fibers are the nodal elliptic curves with Euler
characteristic equal to 1. Hence, there are 12 special fibers.

However, as remarked in Section \ref{sM11}, each point of  $\Mb_{1,1}$
is really a half--point because of  automorphism of order 2 of 
any pointed genus $1$ curve. Therefore, the $24=2\cdot 12$ is the
true degree of the map $\phi$. By the push-pull formula we thus
obtain
$$
\int_{\Mb_{1,1}} \psi_1  = \frac{1}{\deg \phi} 
\int_B \phi^* \, \psi_1  = 
 \frac{1}{\deg \phi} \, (- (S_1,S_1)) = \frac{1}{24} \,.
$$
An interesting corollary of this computation is that 
if $F\to B$ is a smooth family of $1$-pointed genus 1
stable curves over a smooth complete curve $B$
then the set-theoretic degree of the induced map $B\to\Mb_{1,1}$ 
has to be divisible by $12$.

\subsection{}

It is difficult to imagine being able to compute many
intersections of the $\psi$-classes in the above manner.
To begin with, it is essentially impossible to write  
down a sufficiently explicit family of general high 
genus curves, see the discussion in Chapter 6F of 
\cite{HM}. It is therefore amazing that there exist 
several complete and beautiful descriptions of 
the all  possible 
intersection numbers of the form
\begin{equation}
\lan \tau_{k_1} \dots \tau_{k_n} 
\ran \overset{\textup{def}}= 
\int_{\Mb_{g,n}} \psi_1^{k_1} \cdots \psi_n^{k_n} \,, \quad
k_1+\cdots+k_n = 3g-3+n \,.
\label{tau}
\end{equation}
The most striking description was conjectured by Witten \cite{W}
and says the exponential of the
following generating function for the numbers
\eqref{tau} 
\begin{equation}
  \label{free}
    F(t_1,t_2,\dots) = \sum_n \frac{1}{n!} \sum_{k_1,\dots,k_n}
\lan \tau_{k_1} \dots \tau_{k_n} 
\ran \, t_{k_1} \cdots t_{k_n} 
\end{equation}
is a $\tau$-function for the KdV hierarchy of
differential equations. This conjecture was motivated by the 
(physical) analogy with the random matrix models of quantum 
gravity and, in fact, the $\tau$-function thus obtained is the
same as the one that arises in the double scaling of the 1-matrix model
(and discussed in other lectures of this school). The KdV equation 
and the string equation satisfied by the $\tau$-function 
uniquely determine all numbers \eqref{tau}. Alternatively,
the numbers \eqref{tau} are uniquely determined by the 
associated Virasoro constraints. Further discussion can be
found, for example, in \cite{D}. 

\subsection{}
Kontsevich in \cite{K} obtained the KdV equations for 
\eqref{tau} from a combinatorial formula for the 
following (somewhat nonstandard) generating function 
\begin{equation}
K_{g,n}(z_1,\dots,z_n) = \sum_{k_1+\cdots+k_n = 3g-3+n} 
\lan \tau_{k_1} \dots \tau_{k_n} \ran \, \prod \frac{(2k_i-1)!!}{z_i^{2k_i+1}}
\,,\label{KK}
\end{equation}
for the numbers \eqref{tau} with fixed $g$ and $n$. 

\begin{figure}[!hbt]
\centering
\scalebox{.7}{\includegraphics{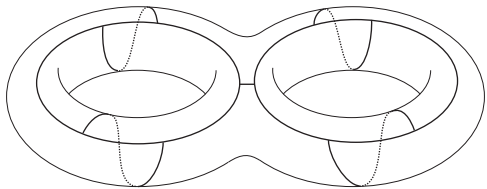}}
\caption{A trivalent map on a genus $2$ surface}
\label{fig11}
\end{figure}

The main ingredient in Kontsevich's combinatorial
formula is a 3-valent graph $G$ embedded in the 
a topological surface $\Sigma_g$. A further 
condition on this graph $G$ is that the complement
$\Sigma_g \setminus G$ is a union of $n$ topological
disks (in particular, this forces $G$ to be connected). 
These disks, called \emph{cells}, have to (bijectively) numbered by
$1,\dots,n$. Two such graphs $G$ and $G'$ are identified
if there exist an orientation preserving homeomorphism of
$\Sigma_g$ that takes $G$ to $G'$ and preserves the labels
of the cells. In particular, every graph $G$ has an  
automorphism group $\Aut G$, which is finite and only seldom
nontrivial. Let $\bG^3_{g,n}$ denote the set of distinct
such graphs $G$ with given values of $g$ and $n$; this is a finite
set. An example of an element of $\bG^3_{2,3}$ is shown in 
Figure \ref{fig11}. 

Another name for a graph $G\subset \Sigma_g$ such that
$\Sigma_g \setminus G$ is a union of cells is a \emph{map}
on $\Sigma_g$. One can imagine that the cells are the 
countries in which the graph $G$ divides the surface $\Sigma_g$.

Kontsevich's combinatorial formula for the function 
\eqref{KK} is the following:
\begin{multline}
  \label{KF}
  K_{g,n}(z_1,\dots,z_n) = \\ 
2^{2g-2+n} \sum_{G\in \bG^3_{g,n}} \frac1{|\Aut G|} 
\, \prod_{\textup{edges $e$ of $G$}} \frac1{z_{\textup{one side of $e$}}+
z_{\textup{other side of $e$}}} \,, 
\end{multline}
where the meaning of the term 
$$
z_{\textup{one side of $e$}}+
z_{\textup{other side of $e$}}
$$
is the following. Each edge $e$ of $G$ separates two cells (which 
may be identical). These cell carry some labels, say, $i$ and $j$.
Then $(z_i+z_j)^{-1}$ is the factor in \eqref{KF} corresponding to  
the edge $e$. 

\begin{figure}[!hbt]
\centering
\scalebox{.7}{\includegraphics{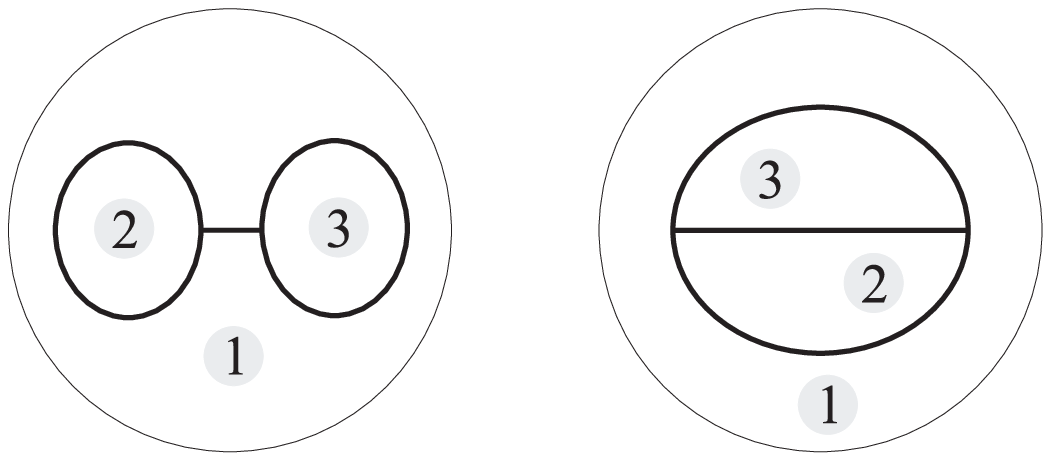}}
\caption{Elements of $\bG^3_{0,3}$}
\label{fig12}
\end{figure}

To get a better feeling for how this works let us look at
the cases $(g,n)=(0,3),(1,1)$ that we understand well. 
The space $\Mb_{0,3}$ is a point and the only nontrivial
integral over it is 
$$
\lan \tau_0 \,  \tau_0\, \tau_0 \ran = \int_{\Mb_{0,3}} 1 = 1 \,.
$$
Thus,
$$
K_{0,3} = \frac1{z_1 z_2 z_3} \,.
$$
The combinatorial side of Kontsevich's formula, however, is
not quite trivial. The set $\bG^3_{0,3}$ consists of 4 elements.
Two of them are shown in Figure \ref{fig12}; 
the other two are
obtained by permuting the cell labels of the graph of the left. 
All these graphs have only trivial automorphisms. Hence, we get
\begin{multline*}
  K_{0,3} = 2 \left( \frac1{2z_1(z_1+z_2)(z_1+z_3)} +
\textup{permutations} \right. \\
+ \left. \frac1{(z_1+z_2)(z_1+z_3)(z_2+z_3)} \right)\,,
\end{multline*}
and, indeed, this simplifies to $(z_1 z_2 z_3)^{-1}$. 
What is apparent in this example is that it is rather
mysterious how \eqref{KF}, which a priori is only a 
rational function of the $z_i$'s, turns out to be 
a polynomial in the variables $z_i^{-1}$. 

\begin{figure}[!hbt]
\centering
\scalebox{.7}{\includegraphics{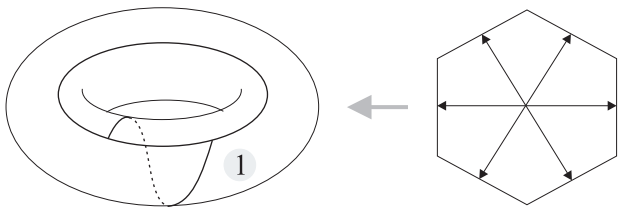}}
\caption{The unique graph in $\bG^3_{1,1}$}
\label{fig13}
\end{figure}

Perhaps this example created a somewhat wrong impression
because in this case \eqref{KF} was much more complicated
than \eqref{KK}. So, let us consider the case $(g,n)=(1,1)$,
where the computation of the unique integral
$$
\lan \tau_1 \ran = \frac{1}{24} 
$$
already does require some work. The unique element of $\bG^3_{1,1}$ 
is shown in Figure \ref{fig13}. 
This graph can be obtained
by gluing the opposite sides of a hexagon, which also 
explains why the automorphism group of this graph is the
cyclic group of order 6 (acting by rotations of the 
hexagon). Thus, \eqref{KF} specializes in this case to 
$$
2 \, \frac{1}{6} \, \frac1{(z_1+z_1)^3}  = \frac{1}{24} \, 
\frac1{z_1^3} \,,
$$ 
as it should.

\subsection{}
Kontsevich was led to the formula \eqref{KF} by 
considering a cellular decomposition of $\M_{g,n}$ 
coming from Strebel differentials. In these 
lectures we shall explain, following \cite{OP}, 
different approach to the formula \eqref{KF} via
the asymptotics in the Hurwitz problem of 
enumerating branched covering of $\PP^1$. This
approach is based on the relation between 
the Hurwitz problem and intersection theory 
on $\Mb_{g,n}$ discovered in \cite{ELSV,FP} and 
on the asymptotic analysis developed in \cite{O1}. 
It has several advantages over the approach 
based on Strebel differentials.

\section{Hurwitz problem}

\subsection{}

Intersection theory on $\Mb_{g,n}$ is about enumerative
geometry of families of stable $n$-pointed curves of 
genus $g$. The significance of the 
space $\Mb_{g,n}$ is that its geometry captures some
essential information about all possible families of 
curves. Through the space $\Mb_{g,n}$, one can learn something about
curves in general from any specific enumerative
problem. If the specific enumerative problem is
sufficiently rich, one can gather a lot
of information about  intersection theory on $\Mb_{g,n}$
from it. Potentially, one can get a complete
understanding of the whole intersection theory, which  
then can be applied to any other enumerative problem. 

Our strategy will be to study such a particular yet
representative enumerative problem. 
 This specific problem will be the Hurwitz problem about 
branched covering of $\PP^1$. That there exists a 
direct connection between Hurwitz problem and 
the intersection theory on $\Mb_{g,n}$
was first realized in \cite{ELSV,FP}. The beautiful formula
of \cite{ELSV} for the Hurwitz numbers will be the basis 
for our computations. 

In fact, we will see that the (exact) knowledge of 
the numbers \eqref{tau} is equivalent to the
\emph{asymptotics} in the Hurwitz problem. This is,
in some sense, very fortunate because asymptotic
enumeration problems often tend to be more
structured and accessible than exact enumeration.

\subsection{}

It is a century-old theme in combinatorics to enumerate
branched coverings of a Riemann surface by another
Riemann surface (an example of which is shown 
schematically in Figure \ref{fig15}). Given degree $d$, positions 
of ramification points downstairs, and their types 
(that is, given the conjugacy class in $S(d)$ of the monodromy 
around each one of them), there exist only finitely many 
possible coverings and the natural question is: how many ?
This very basic enumerative problem arises all over
mathematics, from complex analysis to ergodic theory.
These numbers of branched coverings are directly connected to other 
fundamental objects in combinatorics, namely to
the class algebra of the symmetric group and ---
via the representation theory of finite groups ---
to the characters of symmetric groups. 

We also mention that there is a general, and explicit, correspondence
between enumeration of branched covering of a curve and the 
the Gromov-Witten theory of the same curve, see \cite{OP2}.
From this point of view, the computation of the 
numbers \eqref{tau}, that is, the Gromov-Witten theory of 
a point, arises as a limit in the Gromov-Witten theory of
$\PP^1$ as the degree goes to infinity. This is parallel
to how the free energy \eqref{free} equation arises as the limit in 
the $1$-matrix model.

\subsection{}

\begin{figure}[!hbt]
\centering
\scalebox{.5}[0.25]{\includegraphics{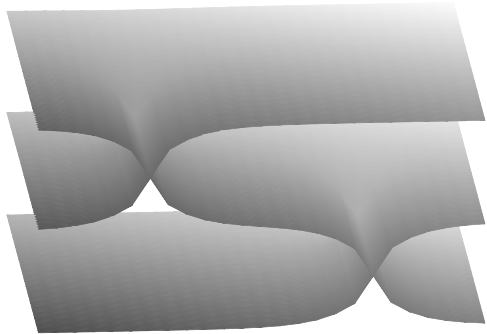}}
\caption{A Hurwitz covering with $\mu=(3)$}
\label{fig14}
\end{figure}

The particular branched covering enumeration problem that 
we will be concerned with can be stated as follows.
The data in the problem are a partition $\mu$ and genus
$g$. Let 
$$
f: C \to \PP^1
$$
be a map of degree
$$
d=|\mu|=\sum \mu_i\,,
$$
where $C$ is smooth connected complex curve of genus $g$.
We require that $\infty\in\PP^1$ is a critical value of 
the map $f$ and the corresponding monodromy has cycle
type $\mu$. Equivalently, this can be phrased as the
requirement that 
divisor $f^{-1}(\infty)$ has the form 
$$
f^{-1}(\infty) = \sum_{i=1}^n \mu_i \, [p_i]\,,
$$
where $n=\ell(\mu)$ is the length of the partition $\mu$ and
$p_1,\dots,p_n\in C$ are the points lying over 
$\infty \in \PP^1$. We further require that 
all other critical values of $f$ are distinct and
nondegenerate. In other words, the map $f^{-1}$ has
only square-root branch points in $\PP^1\setminus\{\infty\}$.
The number $r$ of such square-root branch points is 
given by the Riemann--Hurwitz formula
\begin{equation}
r= 2g-2 + |\mu| + \ell(\mu)  \,.\label{r=}
\end{equation}
An example of such a covering can be seen in Figure 
\ref{fig14} where
$\mu=(3)$ and $r=2$, hence $d=3$, $n=1$, and $g=0$. 

We will call a covering satisfying the above conditions
a \emph{Hurwitz covering}. 
Once the positions of the $r$ simple branchings are
fixed, there are only finitely many Hurwitz coverings provided 
we identify two coverings
$$
f:C\to \PP^1\,, \quad f':C'\to\PP^1
$$
for which there exists an isomorphism $h: C\to C'$ such
that $f=f'\circ h$. Similarly, we define automorphisms
of $f$ as an automorphisms $h: C \to C$ such that $f=f\circ h$.
We will see that, with a very rare exception, Hurwitz coverings
have only trivial automophisms.

By definition, the \emph{Hurwitz number} $\Hur_g(\mu)$ is the number of 
isomorphism classes of Hurwitz coverings with given positions
of branch points. In the special case when such a covering has
a nontrivial automorphism, it should be counted with multiplicity $\frac12$. 

\subsection{}

The Hurwitz problem can be restated as a problem about factoring
permutations into transpositions. This goes as follows. 

Let us pick a point $x\in\PP^1$ which is not a ramification 
point. Then, by basic topology, all information about the 
covering is encoded in the homomorphism
$$
\pi_1(\PP^1\setminus\{\textup{ramification points}\},x) \to
\Aut f^{-1}(x)\cong S(d) \,.
$$
The identification of $\Aut f^{-1}(x)$ with $S(d)$ here is 
not canonical, but it is convenient to pick any one of the $d!$
possible identifications. Then, by construction, the loop around $\infty$ goes
to a permutation $s\in S(d)$ with cycle type $\mu$ and loops around
finite ramification points correspond to some transpositions
$t_1,\dots,t_r$ in $S(d)$. 

The unique relation between those loops in 
$\pi_1$ becomes the equation 
\begin{equation}
t_1 \cdots t_r = s\,.\label{tpr}
\end{equation}
This establishes the equivalence of the
Hurwitz problem with the problem of factoring general permutations 
into transpositions (up to conjugation, since we picked an
arbitrary identification of  $\Aut f^{-1}(x)$ with $S(d)$). 
More precisely, the Hurwitz number $\Hur_g(\mu)$ is the 
number (up to conjugacy, and possibly with an automorphism
factor) of factorization of the form \eqref{tpr} that
correspond to a connected branched covering. A branched
covering is connected when we can get from any point
of $f^{-1}(x)$ to any other point by the action of 
the monodromy group. Thus, the transpositions $t_1,\dots,t_n$
have to generate a transitive subgroup of $S(d)$, which 
is then automatically forced to be the whole of $S(d)$.

The fact that  $t_1,\dots,t_n$ generate $S(d)$ greatly 
constraints the possible automorphisms of $f$. Indeed, the
action of any nontrivial automorphism on $f^{-1}(x)$ has to 
commute with $t_1,\dots,t_n$, and hence with $S(d)$, which is
only possible if $d=2$. 

By the usual inclusion--exclusion principle, it is clear that
one can go back and forth between enumeration of connected and
possibly disconnected coverings. Thus, the Hurwitz problem is 
essentially equivalent to decomposing the powers of one
single element of the 
class algebra of the symmetric group, namely of 
the conjugacy class of a transposition
\begin{equation}
\sum_{1\le i<j\le d} (ij) 
\label{2cycl}
\end{equation}
in the standard conjugacy class basis. There is a classical
formula, going back to Frobenius, for all such expansion
coefficients in terms of irreducible characters. The character
sums that one thus obtains can be viewed as finite analogs
of Hermitian matrix integrals, with the dimension of a
representation $\lambda$ playing the role of the Vandermonde
determinant
and the central character of \eqref{2cycl} in the representation 
$\lambda$ 
playing the role of the Gaussian density, see, for example
\cite{O2,O4} for a further discussion of properties of such sums.

\subsection{}

For us, the crucial property of the Hurwitz problem is its
connection with the intersection theory on the Deligne-Mumford
spaces $\Mgnb$. This connection was discovered, independently,
in \cite{FP} and \cite{ELSV}, the latter paper containing the following
general formula
\begin{equation}
\Hur_g(\mu) = \frac{r!}{|\Aut\mu|} 
\prod_{i=1}^n \frac{\mu_i^{\mu_i}}{\mu_i !} \,
\int_{\Mgnb} \frac{1-\lambda_1+\dots\pm \lambda_g}{\prod (1-\mu_i \psi_i)}\,,
\label{ELSV}
\end{equation}
where $r$ is number of branch points given by \eqref{r=}, $n=\ell(\mu)$
is the length of the partition $\mu$, $\Aut \mu$ is the stabilizer of 
the vector $\mu$ in $S(n)$,
$$
\lambda_i \in H^{2i} \left(\Mgnb\right)\,,\quad i=1,\dots,g\,, 
$$
are the Chern classes of the Hodge bundle over $\Mgnb$ (it is not
important for what follows to know what this is), and finally, the 
denominators are supposed to be expanded into a geometric series
\begin{equation}
\frac{1}{1-\mu_i \psi_i} = 1+ \mu_i \psi_i + \mu_i^2 \psi_i^2 + \dots\,,
\label{gc}
\end{equation}
which terminates because $\psi_i \in H^{2}\left(\Mgnb\right)$
is nilpotent. 

In particular, the integral in the ELSV formula \eqref{ELSV} is
a polynomial in the $\mu_i$'s. The monomials in this polynomial are 
obtained by picking a term in the expansion \eqref{gc} for 
each $i=1,\dots,n$ and then adding a suitable $\lambda$-class
to bring the total degree to the dimension of $\Mgnb$. It is,
therefore, clear that
the top degree term of this polynomial involves only intersections
of the $\psi$-classes and no $\lambda$-classes. That is, 
\begin{equation}
\int_{\Mgnb} = \sum_{k_1+\cdots+k_n=3g-3+n} \prod \mu_i^{k_i} \, 
\lan \tau_{k_1} \dots \tau_{k_n} \ran + \textup{lower degree} \,.\label{topE}
\end{equation}
These top degree terms are precisely the numbers \eqref{tau}
that we want to understand.

\subsection{}\label{sLH}

A natural way to infer something about the top degree part of
a polynomial is to let its arguments go to
infinity. The behavior of the prefactors in \eqref{ELSV}
is given by the Stirling formula
$$
\frac{m^m}{m!} \sim \frac{e^m}{\sqrt{2\pi m}}\,, \quad 
m\to \infty \,.
$$
Let $N$ be a large parameter and let $\mu_i$ depend on $N$ in such a 
way that 
$$
\frac{\mu_i}{N} \to x_i\,, \quad i=1,\dots,n\,, \quad N\to\infty\,,
$$
where $x_1,\dots,x_n$ are finite. We will also additionally
assume that all $\mu_i$'s are distinct and hence $|\Aut\mu|=1$. 
Then by \eqref{topE} and the
Stirling formula, we have the following asymptotics of
the Hurwitz numbers:
\begin{multline}
  \frac{1}{N^{3g-3+n/2}} \, \frac{\Hur_g(\mu)}{e^{|\mu|} r !} \to \\
  \frac1{(2\pi)^{n/2}} \sum_{k_1+\cdots+k_n=3g-3+n} \prod
  \mu_i^{k_i-\frac12} \, \lan \tau_{k_1} \dots \tau_{k_n} \ran =:
  H_g(x)\,. \label{Hg}
\end{multline}
It is convenient to Laplace transform the asymptotics 
$H_g(x)$. Since
$$
\int_0^\infty e^{-sx} \, x^{k-1/2} \, dx = \frac{\Gamma(k+1/2)}{s^{k+1/2}} =
\sqrt{\pi} \, \frac{(2k-1)!!}{2^k\, s^{k+1/2}} \,,
$$
we get 
\begin{equation}
\int_{\R^n_{>0}} e^{-s\cdot x} \, H_g(x) \, dx =
\sum_{k_1+\cdots+k_n=3g-3+n} 
\lan \tau_{k_1} \dots \tau_{k_n} \ran \, \prod 
\frac{(2k_i-1)!!}{(2s_i)^{k_i+1/2}} \,,\label{LHg}
\end{equation}
which up to the following change of variables
$$
z_i = \sqrt{2 s_i} \,, \quad i=1,\dots,n \,,
$$
is precisely the Kontsevich generating function \eqref{KK}
for the numbers \eqref{tau}. 

Thus, we find ourselves in situation which looks rather
comfortable: the generating function that we seek to 
compute is not only related to a specific enumerative
problem but, in fact, it is the Laplace transform of the 
asymptotics in that enumerative problem. People who do 
enumeration know that asymptotics tends to be simpler
than exact enumeration and, usually, the Laplace (or
Fourier) transform of the asymptotics is the most
natural thing to compute. 

This general philosophy is,
of course, only good if we can find a handle on the 
Hurwitz problem. In the following subsection, we will 
discuss a restatement of the Hurwitz problem in terms
of enumeration of certain graphs on genus $g$ surfaces
that we call branching graphs. This description will
turn out to be particularly suitable for our purposes
(which may not be a huge surprise because, after all,
Kontsevich's formula \eqref{KF} is stated in terms of 
graphs on surfaces).

\subsection{}\label{Hbr}

A very classical way to study branched coverings is to 
cut the base into simply-connected pieces. Over each 
of the resulting  regions the covering becomes trivial, that
is, consisting of $d$ disjoint copies of the region 
downstairs, where $d$ is the degree of the covering. 
The structure of the covering is then encoded in the information on 
how those pieces are patched together upstairs. Typically,
this gluing data is presented in the form of a graph,
usually with some additional labels etc. 

There is, obviously,
a considerable flexibility in this approach and some 
choices may lead to much more convenient graph enumeration 
problems than the others. For the Hurwitz problem, we will
follow the strategy from \cite{A}, which goes as follows. 

\begin{figure}[!hbt]
\centering
\scalebox{.7}{\includegraphics{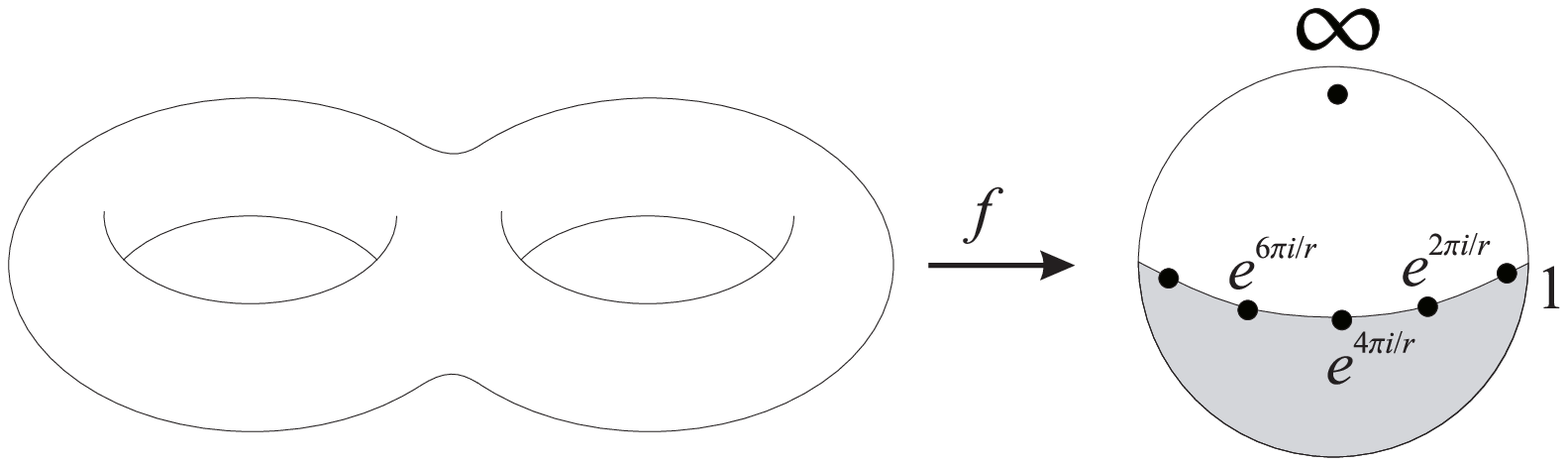}}
\caption{A Hurwitz covering $f:\Sigma_2\to\PP^1$}
\label{fig15}
\end{figure}

Let 
$$
f: C \to \PP^1
$$
be a Hurwitz covering with partition $\mu$ and genus $g$.
In particular, the number $r$ of finite ramification 
points of $f$ is given by the formula \eqref{r=}.  
Without loss of generality, we can assume these 
ramification points to be $r$th roots of unity in $\C$. Let us cut
the base $\PP^1$ along the unit circle $S=\{|z|=1\}$, that is, let
us write
$$
\PP^1 = D_- \sqcup S \sqcup D_+
$$
where 
$$
D_\pm = \{|z|\lessgtr 1\}
$$ 
are the Southern and Northern hemisphere in Figure \ref{fig15},
respectively. 

\begin{figure}[!hbt]
\centering
\scalebox{.5}{\includegraphics{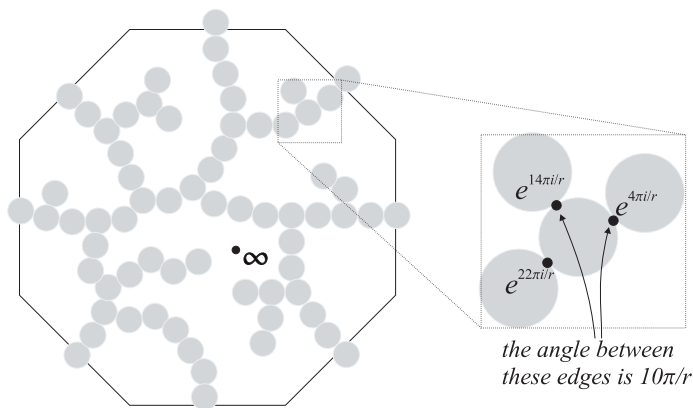}}
\caption{Preimage $f^{-1}(D_-)$  on $\Sigma_2$}
\label{fig16}
\end{figure}

Since the map $f$ is unramified over $D_-$, its preimage 
$f^{-1}(D_-)$ consists of $d$ disjoint disks. Their closures, however,
are not disjoint: they come together precisely at the 
critical points of $f$. By construction, critical points
of $f$ are in bijection with its critical values, that is,
with the $r$th roots of unity in $\PP^1$. Thus, the 
the set $f^{-1}(\overline{D_-}) \subset C$ looks like 
the structure in Figure \ref{fig16}. This structure is, in fact,
a graph $\Gamma$ embedded in a genus $g$ surface. Its vertices
are the components of $f^{-1}(D_-)$ and its edges
are the critical points of $f$ that join those components
together. In addition, the edges of $\Gamma$ (there are $r$ of them)
are labeled by the roots of unity. 

This edge-labeled graph $\Gamma\subset\Sigma_g$ is subject to some 
additional constraints. First, the cyclic order of labels
around any vertex should be in agreement with the 
cyclic order of roots of unity. Next, the complement 
of $\Gamma$ consists of $n$ topological disks, where 
$n$ is the length of the partition $\mu$. Indeed, the 
complement of $\Gamma$ corresponds to $f^{-1}(D_+)$
and $z=\infty$ is the only ramification point in $D_+$. 
The connected components of $f^{-1}(D_+)$ thus
correspond to parts of $\mu$. 

The partition $\mu$ can be reconstructed from the 
edge labels of $\Gamma$ as follows. Pick a cell $U_i$ in 
$f^{-1}(D_+)$. The length of the corresponding part
$\mu_i$ of $\mu$ is precisely the number of times
the map $f$ wraps the boundary $\partial U_i$ around
the circle $S$. As we follow the boundary $\partial U_i$,
we see the edge labels appear in a certain sequence.
As we complete a full circle 
 around $\partial U_i$, the edge labels will make
exactly $\mu_i$ turns around $S$. It is natural to call 
this number $\mu_i$ the \emph{perimeter} of the cell $U_i$.
This perimeter is $(2\pi)^{-1}$ times the sum of \emph{angles} between
pairs of the adjacent edges on $\partial U_i$, where the
angle is the usual angle in $(0,2\pi)$ between the 
corresponding roots of unity, see Figure \ref{fig16}. 

We call a edge-labeled embedded graph $\Gamma$ 
 as above a \emph{branching graph}. By the above
correspondence, the number $\Hur_g(\mu)$ is the 
number of genus $g$ branching graphs with $n$ cells
of perimeter $\mu_1,\dots,\mu_n$. As usual, in the
trivial $d=2$ case, those graphs have to be counted with 
automorphism factors. 

It is this definition of Hurwitz
numbers that we will use for the asymptotic analysis 
in the next lecture.

\subsection{}

\begin{figure}[!hbt]
\centering
\scalebox{.6}{\includegraphics{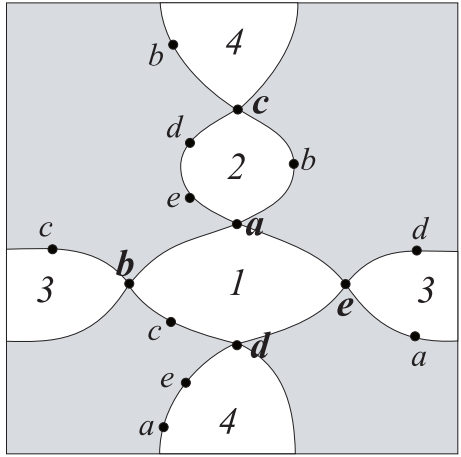}}
\caption{Preimage of $f^{-1}(D_-)$ on the torus $\Sigma_1$}
\label{fig17}
\end{figure}

It may be instructive to consider an example of how this
correspondence between coverings and graphs works. 
Consider the covering corresponding to factorization 
$$
(12)\, (13)\, (24)\, (14)\, (13)=(1243)
$$
of the form \eqref{tpr}. The degree of this covering
is $d=4$, it has $r=5$ ramification points, and the 
monodromy $\mu=(4)$ around infinity. It follows that 
its genus is $g=1$. Let us denote the five finite
ramification points by
$$
\{a,b,c,d,e\}=\{1,e^{2\pi i/5},\dots,e^{8\pi i/5}\} \,.
$$
The preimage of $D_-$ on the torus $\Sigma_1$ consists
of $4$ disks and the momodromies tell us which 
disk is connected to which at which critical point:
for example, at the critical point
lying over $a$, the 1st disk is connected to the 2nd disk. 
This is illustrated in Figure \ref{fig17}
where, among the 3 preimages of any critical value, the one
which is a critical point is typeset in boldface. 
Clearly, any disk in $f^{-1}(D_-)$ has the alphabet $\{a,b,c,d,e\}$
going counterclockwise around its boundary and, in particular,
the cyclic order of the critical values on its boundary
is in agreement with the orientation on $\Sigma_1$.  

\begin{figure}[!hbt]
\centering
\scalebox{.6}{\includegraphics{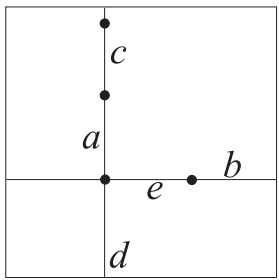}}
\caption{Branching graph corresponding to Figure \ref{fig17}}
\label{fig18}
\end{figure}

Observe that the preimage $f^{-1}(D_+)$ is one cell 
whose boundary
is a 4-fold covering of  the equator. In particular,
the alphabet $\{a,b,c,d,e\}$ is repeated $4$ times around
the boundary of $f^{-1}(D_+)$. Finally, Figure  \ref{fig18}
shows the branching graph 
translation of Figure \ref{fig17}. 

\subsection{}

Finally, a few remarks about how one can prove a formula
like \eqref{ELSV}. This will necessary be a very sketchy
account; the actual details of the proof can be found
in \cite{GV,OP}, as well as in the original paper \cite{ELSV}. 

As mentioned before, the numbers like $\Hur_g(\mu)$ a
special case of in the integrals in the Gromov-Witten 
theory of $\PP^1$, that is, certain intersections on 
the Kontsevich moduli space $\Mb_{g,d,n}(\PP^1)$
of stable degree $d$ maps 
$$
f: C \to \PP^1
$$
from a varying $n$-pointed genus $g$ domain curve $C$ to the 
fixed target curve $\PP^1$. 

Since such a map can
be composed with any automorphism of $\PP^1$, we
have a $\C^\times$-action on $\Mb_{g,d,n}(\PP^1)$. A theory
due to Graber and Pandharipande \cite{GP} explains how
to localize the integrals in Gromov-Witten theory
to the fixed points of the action of the torus $\C^\times$. 
These fixed point loci in $\Mb_{g,d,n}(\PP^1)$ are,
essentially, products of Deligne-Mumford spaces $\Mb_{g_i,n_i}$
for some  $g_i$'s and $n_i$'s. Indeed, 
only very few maps are fixed by the action of the 
torus. Namely, for the standard $\C^\times$-action on $\PP^1$
and an irreducible domain curve $C$ the only choices
are the degree $0$ constant maps 
to $\{0,\infty\}=\left(\PP^1\right)^{\C^\times}$
or the degree $d$ map
$$
\PP^1 \owns z \mapsto z^d \in \PP^1 \,.
$$
In general, the domain curve is allowed to be reducible,
but still any torus-invariant map has to be of the 
above form on each component $C_i$ of $C$. Once all discrete
invariants of the curve $C$ are fixed (that is, the 
combinatorics of its irreducible components, their 
genera and numbers of marked points on them) the 
remaining moduli parameters are only a choice of
a bunch of curves to collapse plus a choice of where 
to attach the non-collapsed $\PP^1$'s to them. That is,
the torus-fixed loci are products of Deligne-Mumford
spaces, modulo possible automorphisms of the combinatorial
structure. 

In this way integrals in the Gromov-Witten 
theory of $\PP^1$ can be reduced, at least in principle,
to computing intersections on $\Mgnb$. 
An elegant localization analysis leading to the 
ELSV formula is presented in \cite{GV}, see also \cite{OP}.

\section{Asymptotics in Hurwitz problem}

\subsection{}

Our goal now is to see how the Laplace transform \eqref{LHg}
of the asymptotics \eqref{Hg} in the Hurwitz problem 
turns into Kontsevich's combinatorial formula \eqref{KF}.
The formulation of the Hurwitz problem in terms of 
branching graphs, see Section \ref{Hbr}, looks promising.
Indeed, a branching graph $\Gamma$ is by definition
 embedded in a topological genus $g$ surface $\Sigma_g$ and
and it cuts $\Sigma_g$ into $n$ cells. Here the numbers $g$ and $n$ are
the same as the indices in $\Mgnb$, on the intersection theory on which
we are trying to understand. Similarly, in Kontsevich's
formula we have a graph $G$ embedded into $\Sigma_g$ and
cutting it into $n$ cells. This graph $G$, however, is a more
modest object: it does not have any edge labels and it
is allowed to have only $3$-valent vertices. 

Recall that we denote by $\bG^3_{g,n}$ the set of all possible
$3$-valent graphs as in Kontsevich's formula \eqref{KF}. 
Let us introduce two larger sets
$$
\bG^3_{g,n} \subset \bG^{\ge 3}_{g,n} \subset \bG_{g,n} \,,
$$ 
on which, by definition, the $3$-valence condition is 
weakened to allow vertices of valence $3$ or more, and 
dropped altogether, respectively. The elements of $\Ggnb$ can be obtained from
elements of $\bG^3_{g,n}$ by contracting some edges. In particular,
the set $\Ggnb$ is still a finite set. Similarly, denote
by $\bH_{g,\mu}$ the set of all branching graphs with given
genus $g$ and perimeter partition $\mu$. Our first order
of business is to construct a map 
$$
\bH_{g,\mu}  \to \Ggnb \,,
$$ 
which we call the \emph{homotopy type} map. This map
is the composition of the map 
$$
 \bH_{g,\mu}  \to \bG_{g,n}
$$
which simply forgets the edge labels with the 
map
$$
\bG_{g,n} \to \Ggnb \,,
$$
which does the following. First, we remove all 
univalent vertices together with the incident edge. 
After that, we remove the all 
remaining $2$-valent vertices joining their two
incident edges. What is left, by construction, 
has only vertices of valence $3$ and higher and
still cuts $\Sigma_g$ into $n$ cells. 

\begin{figure}[!hbt]
\centering
\scalebox{.5}{\includegraphics{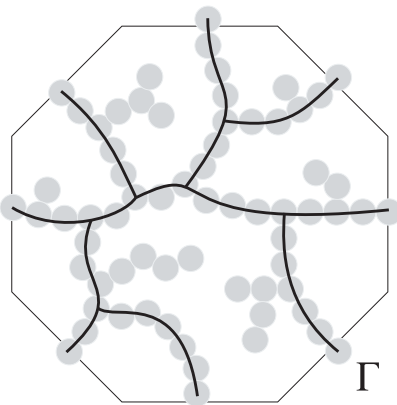}}
\caption{The homotopy type of the branching graph  from Figure \ref{fig16}}
\label{fig19}
\end{figure}

We remark that in the two exceptional cases
$(g,n)=(0,1),(0,2)$, which correspond to unstable
moduli spaces, what we get in the end (a point and
a circle, respectively) is
not really an element of $\Ggnb$. In all other 
cases, however, we do get an honest element of $\Ggnb$. 
Figure \ref{fig19} illustrates this procedure
applied to the branching graph from Figure \ref{fig16}. 

\subsection{}

Now let us make the following simple but 
important observation. Since the set $\Ggnb$ is 
\emph{finite} and we are interested in the asymptotics
of $\Hur_{g}(\mu)$ as $\mu\to\infty$ while
keeping $g$ and $n$ fixed, we can just do the
asymptotics separately for each homotopy type
and then sum over all possible homotopy types. 
The Laplace transform \eqref{LHg} will then
be also expressed as a sum over all corresponding
homotopy types $G$ in $\Ggnb$.

We now claim that not only Kontsevich's combinatorial formula
\eqref{KF} is the Laplace transformed asymptotics
\eqref{LHg} but, in fact, the summation over $G\in \bG^3_{g,n}$
in Kontsevich's formula corresponds precisely to 
summation over possible homotopy types. Since 
there are non-trivalent homotopy types, implicit in 
this claim is the statement that \emph{non-trivalent homotopy types
do not contribute to asymptotics}. 

\subsection{}

What do we need to do to get the asymptotics of the 
number of branching graphs of a given homotopy type $G$ ?
What would suffice is to have a simple way to 
enumerate all such branching graphs. To enumerate
all branching graphs with given homotopy type $G$, we
need to retrace the steps of the homotopy type map. 
Imagine that the homotopy type graph $G$ is a fossil
from which we want to reconstruct some prehistoric
branching graph $\Gamma$. What are the all possible
ways to do it ?

The answer to all these rhetoric questions is quite
simple. It is easy to see that the preimage of 
any edge in $G$ is some subtree in the original
branching graph $\Gamma$.  In addition, all
these trees carry edge-labels which were erased
by the homotopy type map. Thus, for any edge
$e$ of $G$, we need to take a tree $T_e$ whose edges
are labeled by roots of unity.  In particular, there is a 
canonical way to make this tree planar, that is,
embed it in the plane in such a way that the
cyclic order of edges around each vertex 
agrees with the order of their labels. 
In particular, each such tree is a \emph{branching
tree}, that is, it satisfies the $(g,n)=(0,1)$ case
of our definition of a branching graph\footnote{
A small and inessential detail is that the labels $T_e$ are 
taken from a larger set of roots of unity}.

\begin{figure}[!hbt]
\centering
\scalebox{.6}{\includegraphics{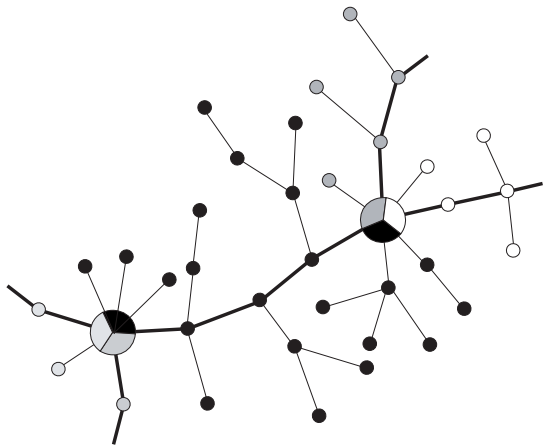}}
\caption{Assembling a branching graph from edge trees}
\label{fig20}
\end{figure}

Next, these trees are to be glued into 
the graph $\Gamma$ by identifying some of their vertices, as in 
Figure \ref{fig20}. This means that each of these
branching trees carries two special vertices, which we call 
its \emph{root} and \emph{top}. These special
vertices of $T_e$  mark the
places where $T_e$ is attached to the other trees in $\Gamma$.
We will call such a branching tree with two marked
vertices an \emph{edge tree}.

\subsection{}

Now we have a procedure which from a homotopy type $G$
and a collection of edge trees $\{T_e\}$ with distinct
labels assembles a branching graph $\Gamma$. This 
procedure, which we will call \emph{assembly},
does have some imperfections. Those imperfections
will be discussed momentarily, but first we want
to make the following important observation.

Since the homotopy type graph $G$ is something fixed and 
finite, \emph{the whole asymptotics of the branching 
graphs lies in the edge trees}. For a large random 
branching graph $\Gamma$, those edge trees will
be large random trees. This is how the theory
of random trees enters the scene. Fortunately,
a large random tree is a very well studied and a very  
nicely behaved object, see for example \cite{Pit} for a
particularly enjoyable introduction. It turns out that
all the information we need about random trees is 
either already classical or can be easily deduced from 
known results. 

In fact, all required knowledge about random trees  
can be quite easily deduced (as was done in 
\cite{OP}) from the first principles, which in this case, is the
following formula going back to Cayley \cite{Sta}. Consider all 
possible trees $T$ with the vertex set $\{1,\dots,m\}$.
For any such tree $T$, we have a function $\val_T(i)$
which takes the vertex $i=1,\dots,m$ to its valence in $T$. 
The information about all vertex valences in all possible
trees $T$ is encoded in the following generating
function 
\begin{equation}
\sum_T z_1^{\val_T(1)} \cdots z_m^{\val_T(m)} =
z_1 \cdots z_m (z_1 + \cdots + z_m)^{m-2} \,.
\label{Cay}
\end{equation}
A probabilistic restatement of this result is
the following. The valence $\val_T(i)$ is the
number of edges of $T$ incident to the 
vertex $i$. Let us cut all edges in half;
since there were $m-1$ edges of $T$, we 
get $2m-2$ half edges. The formula \eqref{Cay}
says that the same distribution of half-edges
can be obtained as follows: give every 
vertex a half-edge and the remaining $m-2$ 
edges just throw at the 
vertices randomly like darts. 

What is then the valence of a given vertex in 
a random tree $T$ ? It is 1 for the half edge
allowance that it always gets plus its share in the
random distribution of $m-2$ darts among 
$m$ targets. As $m\to\infty$, this share goes to 
a Poisson random variable with mean $1$. 
In other words, as $m\to\infty$ we have
\begin{equation}
\Prob\{\val_T(i)=v\} \to \frac{e^{-1}}{(v-1)!} \,,
\quad v=1,2,\dots \,.\label{valP}
\end{equation}
For different vertices,
their valences become independent in the
$m\to\infty$ limit. 

Also, setting all variables in \eqref{Cay} to $1$ we find
that the total number of trees with 
vertex set $\{1,\dots,m\}$ is $m^{m-2}$.

\subsection{}\label{assc}

Now it is time to talk about how the assembly map 
differs from being one-to-one (it clear that it is 
onto). 

First, it may happen that the cyclic order of edge labels
is violated at one of the vertices of $G$ where we 
patch together different edge trees. If this is the case, we simply
declare the assembly to be a failure and do nothing. 
The probability of such an assembly failure in the 
large graph limit can be computed as follows. 
Suppose that we need to glue together three vertices
with valences $v_1$, $v_2$, and $v_3$. From \eqref{valP}, the 
chance of seeing these particular valences is 
$$
\frac{e^{-3}}{(v_1-1)! (v_2-1)! (v_3-1)!} \,.
$$
On the other hand, the conditional
probability that the edge labels in the resulting
graph are cyclically ordered, given that they were
cyclically ordered before gluing is easily
seen to be
$$
\frac{(v_1-1)! (v_2-1)! (v_3-1)!}{(v_1+v_2+v_3-1)!} \,.
$$
Hence the success rate of the assembly at a particular
trivalent vertex is 
$$
e^{-3}\sum_{v_1,v_2,v_3\ge 1} \frac1{(v_1+v_2+v_3-1)!} =
\frac{e^{-2}}{2} \,.
$$
Assembly failures at distinct vertices being asymptotically
independent events, this goes into an overall factor and,
eventually in the prefactor in \eqref{KF}.

At this point it should be clear that there in no 
need to consider nontrivalent vertices. Indeed, a homotopy
type graph with a vertex of valence $\ge 4$ can be 
obtained from a trivalent graph by contracting some 
edges, hence corresponds to the case when some
of the edge graphs are trivial. It is obvious that
the chances that a large random tree came out 
empty are negligible. Hence, nontrivalent graphs make
indeed no contribution to the asymptotics and can
safely be ignored.

\subsection{}

The second (minor) issue with the assembly map is 
that we can get the same, that is, isomorphic
branching graphs starting from different collections
of the edge trees. This happens if the homotopy
type graph $G$ has nontrivial automorphisms.
It is clear that the group $\Aut(G)$ acts on 
edges of $G$ and, hence, acts by permutations on 
collections of edge trees preserving the 
isomorphism class of the assembly output. 
It is also clear that  the 
chance for a large edge tree to be isomorphic to 
another edge tree (or to itself with root and
top permuted) is, asymptotically, zero. Hence
almost surely this $\Aut(G)$ action is free and
hence there is an overcounting of branching
graphs by exactly a factor of $|\Aut(G)|$. 
This explains the division by $|\Aut(G)|$ in 
\eqref{KF}.

\subsection{}

Now, after explaining the summation over
trivalent graphs and the automorphism 
factor in \eqref{KF}, we get to the heart
of Kontsevich's formula --- the product 
over the edges. 

It is at this point that the convenience
(promised in Section \ref{sLH})
of working with the Laplace transform \eqref{LHg}
rather then the asymptotics \eqref{Hg} itself
can be appreciated.
We will see shortly that, asymptotically,
the cell perimeters of a branching graph
$\Gamma$ assembled from a 3-valent
graph $G$ and bunch of 
random edge trees $\{T_e\}$ is a sum of 
independent contributions from each edge of $G$. 
This makes the Laplace transform \eqref{LHg}
factor over the edges of $G$ as in \eqref{KF}.
To justify the above claim, we need to take a closer
look at a large typical edge tree.

\subsection{}

\begin{figure}[!hbt]
\centering
\scalebox{.6}{\includegraphics{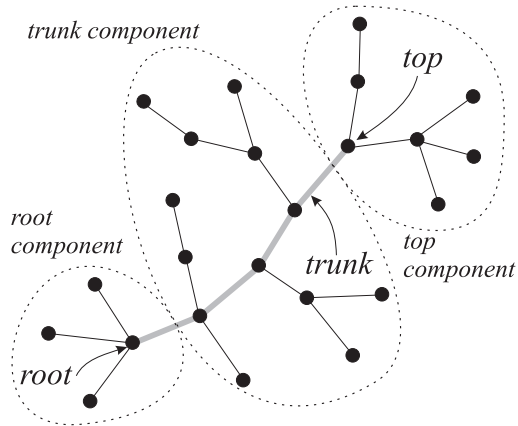}}
\caption{The components of an edge tree}
\label{fig21}
\end{figure}

Let $T$ be an edge tree. 
It has two marked vertices, root and top;
let us call the path joining them the 
\emph{trunk} of $T$. The tree $T$ naturally 
splits into 3 parts: 
the root component, the top component,  and the trunk component,
according to their closest trunk point. 
This is illustrated in Figure \ref{fig21}. 

Figure \ref{fig21} may give a wrong idea of the relative
size of these components for a typical large edge 
tree. Let $M\to\infty$ be the size (e.g.\  the number
of vertices) in $T$. It is known and can be 
without difficulty deduced from \eqref{Cay} (see 
for example \cite{OP}) that the size $M$ distributes itself 
among the three components of $T$ in the $M\to\infty$ limit
as follows.

First, the size of the root and the top component 
stays finite in the $M\to\infty$ limit. In fact,
it goes to the \emph{Borel distribution}, given
by the following formula
$$
\Prob(k) =  \frac{k^{k-1}\, e^{-k}}{k!}\,, 
\quad k=1,2,\dots \,.
$$
Second, the typical size of the trunk is 
of order $\sqrt{M}$. More precisely, scaled 
by  $\sqrt{M}$, the trunk size distribution
goes to the Rayleigh distribution with 
density
$$
x\, e^{-x^2/2} \, dx \,, \quad x\in (0,\infty) \,.
$$
For our purposes, however, it only matters that 
the size of all these parts is $o(M)$ as 
$M\to\infty$. 

The overwhelming majority of vertices lie,
therefore, somewhere in the branches of the 
trunk component of $T$. What is very 
important is that, after assembly, any
such vertex will find itself completely surrounded
by a unique cell.  As a result, it will
contribute exactly $1$ to that cell's perimeter. 
What this analysis shows it that, asymptotically,
the cell perimeters are determined simply by the
the number of such interior trunk vertices
ending up in a given cell, all other contributions
to perimeters being $o(M)$. It should be clear
that such contributions of distinct edges of $G$
are indeed independent, leading to the factorization 
in \eqref{KF}.

\subsection{}

What remains is to determine is what the edge factors
are, that is, to determine the actual 
contribution of an edge tree $T$ to the perimeters
of the adjacent cells. 

All we need to know for this is to know how vertices in the trunk component
distribute themselves between the two sides of the trunk,
as in Figure \ref{fig21}. One shows, see \cite{OP} and below, 
that the fraction of the
vertices that land on a given side of the trunk
is, asymptotically, \emph{uniformly distributed} on $[0,1]$.
This reduces the computation of the edge factor
to computing one single integral. That computation 
will be presented in a moment, after we review 
the knowledge that we have accumulated so far.

\subsection{}

Let $G$ be a 3-valent map with $n$ cells. It has
$$
 |E(G)|= 6g-6+3n 
$$
edges and 
$$
|V(G)|= 4g-4+2n
$$
vertices, which follows from the Euler 
characteristic equation
$$
|V(G)| - |E(G)| + n = 2- 2g
$$
combined with the 3-valence condition $3|V(G)|=2|E(G)|$. 

Let $e\in E(G)$ be an edge of $G$ and let $T_e$ be
the corresponding edge tree. Let $d_e$ be the number
of vertices of $T_e$. Ignoring the few vertices 
on the trunk itself, the vertices of $T_e$ distribute
themselves between the two sides of the trunk of $T_e$.
Let's say that $p_e$ vertices are on the one side
and define the number $q_e$ by 
\begin{equation}
p_e + q_e =  d_e \,.\label{split}
\end{equation}
It is clear that $q_e$ is the approximate number of 
vertices on the other side of the trunk. 
We call the numbers $p_e$ and $q_e$ the 
\emph{semiperimeters} of the tree $T_e$. 

The basic question, which we now can answer in the 
large graph limit, is
how many branching graphs $\Gamma$ have
given semiperimeters $\{(p_e,q_e)\}_{e\in E(G)}$. 
This distribution can be computed asymptotically as
follows.

\subsection{}

First, there are some overall factors that come 
from automorphisms of $G$ and the assembly 
success rate. Recall that in Section \ref{assc}
we saw that the assembly success rate is 
$e^{-2|V(G)|} 2^{-|V(G)|}$. 

Second, for every edge $e\in E(G)$ we need to 
pick an edge tree $T_e$ with $d_e$ vertices. 
As we already learned from \eqref{Cay}, the number
of vertex-labeled trees with $d_e$ vertices is 
$d_e^{d_e-2}$. Vertex
labels can be traded for edge labels at the expense
of the factor $d_e!/(d_e-1)!=d_e$, hence there are $d_e^{d_e-3}$
edge labeled trees with $d_e$ vertices. 
The choice of the root vertex brings in additional
factor of $d_e$ choices. Once the root is fixed,
the condition \eqref{split} dictates the position of
the top, so there is no additional freedom in 
choosing it. To summarize, there are $\sim d_e^{d_e-2}$ edge
trees with given semiperimeters $p_e$ and $q_e$. 

Third, the edge labels of $\Gamma$ is a shuffle of
edge labels of the trees $T_e$. Let 
$$
r=\sum_{e}(d_e-1)
$$
be the total number of edges in $\Gamma$ (and, hence, also the
total number of simple branch points in the Hurwitz 
covering  corresponding to $\Gamma$). Obviously, there
are 
\begin{equation}
\frac{r!}{\prod_{e\in E(G)} (d_e-1)!} \label{de!}
\end{equation}
ways to shuffle edge labels of $\{T_e\}$ into edge
labels of $\Gamma$. 

Putting it all together, we obtain the following 
approximate expression for the number of branching
graphs with given semiperimeters  $\{(p_e,q_e)\}_{e\in E(G)}$
\begin{equation}\label{r!}
 \frac{r! \, e^{-2|V(G)|} \, 2^{-|V(G)|}}{|\Aut(G)|} \,
\prod_{e\in E(G)} \frac{d_e^{d_e-2}}{(d_e-1) !} \sim 
\frac{r! \, e^d \, 2^{-|V(G)|}}{|\Aut(G)|} \,
\prod_{e\in E(G)} \frac{1}{\sqrt{2\pi}\, d_e^{3/2}} \,,
\end{equation}
where 
$$
d = |V(\Gamma)|=  \sum_{e} d_e - 2 |V(G)|
$$
is the degree of the corresponding Hurwitz covering and
the RHS of \eqref{r!} is obtained from the LHS by the 
Stirling formula. 

Note that the factor $r! \, e^d$ precisely cancels with
prefactor in \eqref{Hg}. 

\subsection{}

Since the cell perimeters of $\Gamma$ are the sums 
of edge tree semiperimeters along the boundaries of 
the cells, the computation of the Laplace transform 
\eqref{LHg} indeed boils
down to the computation of a singe edge factor
\begin{align*}
  \frac{1}{\sqrt{2\pi}} \iint_{p,q>0} \, \frac{e^{-p s_1 - q s_2
      }}{(p+q)^{3/2}} \, dp \, dq & =   \frac{1}{\sqrt{2\pi}} \,
\frac1{s_1-s_2}\int_{x>0} 
\left(e^{-s_1 x} - e^{-s_2 x}\right) \, \frac{dx}{x^{3/2}} \\
&= \frac{1}{\sqrt{2\pi}}\,
\Gamma\left(-\tfrac12\right) \, \frac{\sqrt{s_1} - \sqrt{s_2}}{s_1-s_2} \\
&=\frac{\sqrt{2}}{\sqrt{s_1} +
    \sqrt{s_2}} \,,
\end{align*}
where we set $x=p+q$. 

Recall that the relation between the Laplace transform
variables $s_i$ in \eqref{LHg}  and the variables $z_i$ in Kontsevich's
generation function \ref{KK} is 
$$
z_i = \sqrt{2 s_i} \,, \quad i=1,\dots,n \,. 
$$
Thus we get indeed the LHS of \eqref{KF}, including the
correct exponent of $2$, which is 
$$
|E(G)|- |V(G)| = 2g-2 + n  \,.
$$
This completes the proof of Kontsevich's formula \eqref{KF}.

\section{Remarks}

\subsection{}

Since random matrices are the common thread of many talks at this 
school, let us point out various connections between moduli of 
curves and random matrices. As we already discussed, the original
KdV conjecture of Witten was based on physical parallelism between
intersection theory on $\Mgnb$ and the double scaling limit of 
the Hermitian $1$-matrix model. Despite many spectacular achievements
by physicists as well as mathematicians, this double scaling
seems to remain a source of serious mathematical challenges,
in particular, it appears that no direct mathematical connection 
between it and moduli of curves is known. On the other hand, 
there is a very direct connection between what we did and 
another, much simpler, matrix model, namely, the edge scaling
of the standard GUE model. This connection goes as follows.

Recall that by Wick formula the coefficients of the $1/N$ expansion of the 
following $N\times N$ Hermitian matrix integral 
\begin{equation}
\int e^{-\tr H^2} \, \prod_{1}^{m} \tr H^{k_i} \, dH \, \label{dH}
\end{equation}
are the numbers of 
ways to glue a surface of given topology from $m$ polygons
with perimeters $k_1,\dots, k_n$. The double scaling limit
of the $1$-matrix model is concerned with gluing a given
surface out of a very large number of small pieces. 
An opposite asymptotic regime is when the number $m$ of 
pieces stays fixed while their sizes $k_i$ go to 
infinity. Since for large $k$ the trace $\tr H^k$ picks out
the maximal eigenvalues of $H$, this asymptotic regime
is about largest eigenvalues of a Hermitian random 
matrix. In the large $N$ limit, the distribution of 
largest eigenvalues of $H$ is well known to be the Airy
ensemble. This edge scaling random matrix ensemble is
very rich, yet susceptible to a very detailed mathematical
analysis. In particular, the individual distributions 
of eigenvalues were found by Tracy and Widom in \cite{TW}.
They are given in terms of certain solutions of the 
Painlev\'e II equation. 

The connection between GUE edge scaling and what we were
doing is the following. If one takes a branching graph 
as in Figure \ref{fig16} and strips it off its edge labels,
one gets a map on genus $\Sigma_g$ with a few
cells of very large perimeter, that is, an object
of precisely the kind enumerated by \eqref{dH} in the 
edge scaling regime. We argued that almost all vertices
of a large branching $\Gamma$ graph are completely surrounded
by a unique cell, hence contribute exactly $1$ to that
cell's perimeter regardless of the edge labels. This
shows that edge labels play no essential role in 
the asymptotics, thus establishing a direct connection
between Hurwitz numbers asymptotics and GUE edge
scaling. A similar direct connection can be established
in other situations, for example, between GUE edge
scaling and distribution of long increasing subsequences
in a random permutation, see \cite{O1}. Since a great deal
is known about GUE edge scaling, one can profit 
very easily from having a direct connection to it.
In particular, one can give closed error-function-type
formulas for a natural generating functions (known as
$n$-point functions) for the numbers \eqref{tau}, see \cite{O3}. 

There exists another matrix model, namely the 
Kontsevich's matrix model \cite{K}, specifically
designed to reproduce the graph summation in \eqref{KF}
as its diagrammatic expansion. Once the combinatorial formula 
\eqref{KF} is established, this Kontsevich's model can 
be used to analyze it, in particular, to prove the 
KdV equations, see \cite{K} and also \cite{D}.  

Alternatively, the KdV equations can be pulled back from 
the GUE edge scaling  (where they have been studied
in depth by Adler, Shiota, and van Moerbeke) via the 
above described connection, see the exposition in \cite{O3}.

\subsection{}

In our approach, the intersections \eqref{tau}, the 
combinatorial formula \eqref{KF}, the KdV equations
etc.\ appear through the asymptotic analysis of 
the Hurwitz problem. The ELSV formula \eqref{ELSV},
which is the bridge between enumeration of branched
coverings and the intersection theory of $\Mgnb$ is,
on the hand, an exact formula. It is, therefore,
natural to ask for more exact  bridges between
intersection theory, combinatorics, and 
integrable systems. 

After the moduli space $\Mgnb$ of stable curves,
a natural next step is the Gromov-Witten theory of 
$\PP^1$, that is, the intersection 
theory on the moduli space $\Mgnb(\PP^1,d)$ of stable 
degree $d$ maps
$$
C \to \PP^1
$$
from an $n$-pointed  genus $g$ curve $C$ to the projective line $\PP^1$.
More generally, one can replace $\PP^1$ by some
higher genus target curve $X$. It turns out, see \cite{OP2}, that 
there is a simple dictionary, which we call the 
Gromov-Witten-Hurwitz correspondence, between enumeration of 
branched coverings of $\PP^1$ and Gromov-Witten 
theory of $\PP^1$. This correspondence naturally connects
with some very beautiful combinatorics and 
integrable systems, the role of random matrices now being
played by random partitions. The connection with the integrable
systems is seen best in the equivariant Gromov-Witten 
theory of $\PP^1$, where the 2-Toda lattice hierarchy
of Ueno and Takasaki plays the role that KdV played for
$\Mgnb$.

\end{document}